\documentclass[10pt]{article}  
\usepackage{amsmath}
\usepackage{amssymb}
\usepackage{epsfig}
\usepackage{graphics}
\usepackage{theorem}

\theoremheaderfont{\scshape}

\theorembodyfont{\normalfont}

\numberwithin{equation}{section}

\title{An explanatory model to validate the way water activity rules
periodic terrace generation in \emph{Proteus mirabilis} swarm}
\author{
Emmanuel Fr\'enod \thanks{Universit\'e Europ\'enne de Bretagne, Lab-STICC (UMR CNRS 3192), 
Universit\'e de Bretagne-Sud, Centre Yves Coppens, Campus de Tohannic,
F-56017, Vannes} \and 
Olivier Sire \thanks{Universit\'e Europ\'eenne de Bretagne, LIMATB, 
Universit\'e de Bretagne-Sud, Centre Yves Coppens, Campus de Tohannic,
F-56017, Vannes}
}

\newcommand{\nit}{\mathbb{N}}
\newcommand{\rit}{\mathbb{R}}
\newcommand{\ritd}{{\mathbb{R}^2}}

\newcommand{\tc}{{t}}
\newcommand{\xc}{{x}}
\newcommand{\posmax}{{x_\text{max}}}
\newcommand{\age}{{a}}
\newcommand{\bge}{{b}}

\newcommand{\qc}{{Q}}
\newcommand{\qcinit}{{Q}_0}
\newcommand{\qcz}{\overline{Q}}
\newcommand{\zec}{{\zeta}}
\newcommand{\rhc}{{\rho}}
\newcommand{\ggc}{{G}}
\newcommand{\hc}{{h}}
\newcommand{\hhc}{{H}}
\newcommand{\biom}{{\cal E}}
\newcommand{\bioz}{{\overline{\cal E}}}
\newcommand{\biomm}{{\cal M}}
\newcommand{\biomn}{{\cal N}}
\newcommand{\gamt}{{\gamma_t}}
\newcommand{\gamd}{{\gamma_d}}
\newcommand{\etac}{{\eta}}
\newcommand{\kac}{{\kappa}}
\newcommand{\alc}{{\alpha}}
\newcommand{\muc}{{\mu}}
\newcommand{\xic}{{\xi}}
\newcommand{\tac}{{\tau}}
\newcommand{\chc}{{\chi}}
\newcommand{\trnc}{{\cal T}}
\newcommand{\rrc}{{\cal D}}

\newcommand{\prptab}{{\nu}}

\newcommand{\timest}{{\Delta t}}
\newcommand{\agest}{{\Delta a}}
\newcommand{\bgest}{{\Delta b}}
\newcommand{\posst}{{\Delta x}}

\newcommand{\ds}{\displaystyle}

\newcommand{\fracp}[2]{\frac{\partial #1}{\partial #2}}

\begin{document}
\maketitle

{\small {\bf Abstract - } This paper explains the biophysical principles
which, according to us, govern the \emph{Proteus mirabilis} swarm phenomenon.
Then,  this explanation is translated into a mathematical model, essentially based
on partial differential equations. This model is then implemented using numerical
methods of the finite volume type in order to make simulations. The simulations 
show most of the characteristics which are observed in \emph{situ} and in particular
the terrace generation.
}

{\small {\bf Keywords - } 
\emph{Proteus mirabilis} swarm, modelling, partial differential equations,
finite volumes.
}

\tableofcontents

\section{Introduction}
\label{intr}
\emph{Proteus mirabilis} is a pathogenic  bacterium of the urinary tract. 
Under specific conditions, within  a \emph{Proteus mirabilis} colony
grown on
a solid substratum,
some bacteria, the standard form of which is a short cell,
called ``swimmer'' or ``vegetative'', undergo differentiation into
elongated cells, called ``swarmer'' cells, capable of translocation. Having this 
translocation ability, those swarmer cells may begin to colonize the 
surrounding medium.
After a given period of time, the swarm phenomenon stops and the swarmer cells 
de-differentiate to produce vegetative cells which would, 
later on, be able to differentiate again and resume colonization.

\:

\:

\begin{figure}
\centering 
\includegraphics[height=10cm]{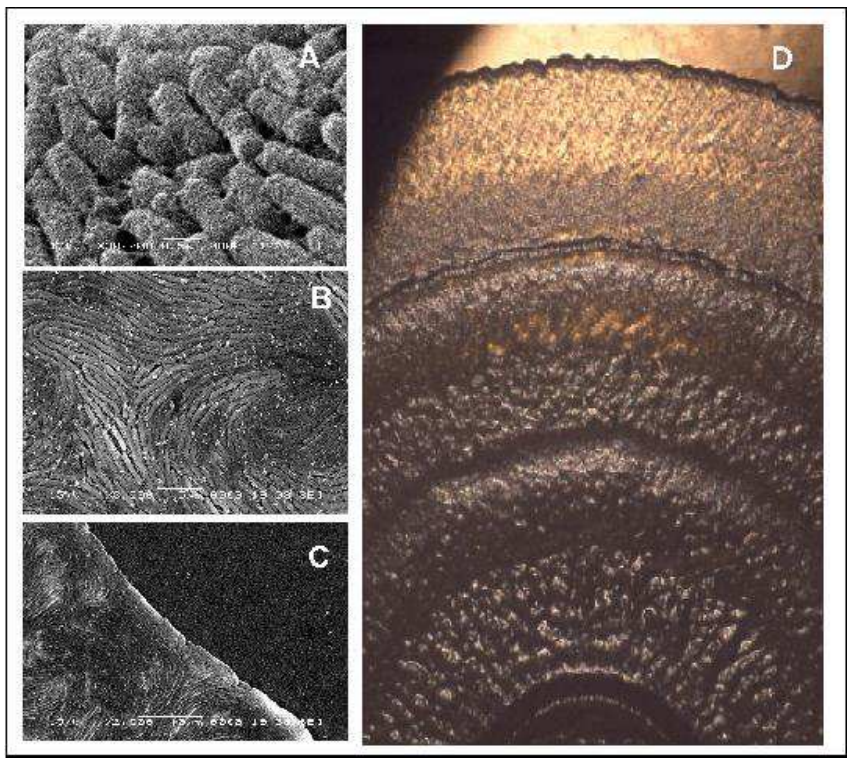}
\caption{\emph{Proteus mirabilis}: Short vegetative cells (A); 
swarmer cells (B); periphery of a colony (C) and a swarm colony with terraces (D).
Scales : (A): bar = 0.5 $\mu m$; (B): bar = 5 $\mu m$; (C): bar = 10 $\mu m$ and 
(D): colony radius = 3 $cm$. }
\label{fig2}
\end{figure}

\begin{figure}
\centering 
~\hspace{-2cm}
\includegraphics[width=18cm]{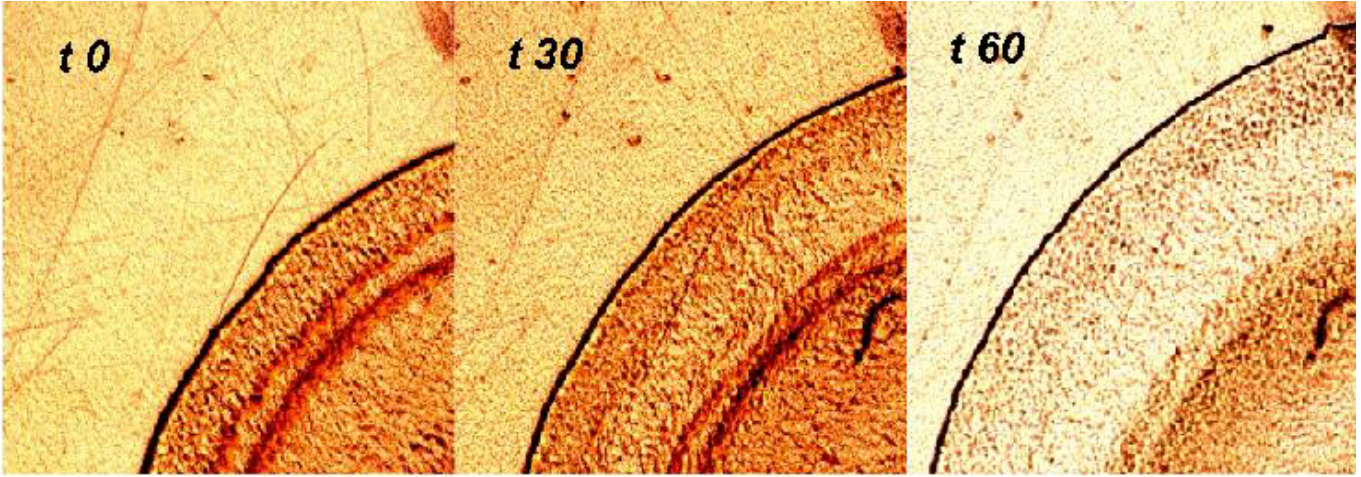}
\caption{Active swarming phase on agar. The figure shows three snapshots 
30 min apart during a migration phase on agar.}\label{figilsw}
\end{figure}

This swarm phenomenon is illustrated in Figures \ref{fig2} and \ref{figilsw}.
In Figure \ref{fig2}, the \emph{Proteus mirabilis} cells can be seen in their
standard vegetative forms (A).
They are about 2 $\mu m$ long. 
The elongated swarmer cells which range between 10 to 25 $\mu m$ 
and which are covered by hundreds of flagella can also
be seen (B). 
The third panel (C) 
displays a close up of the peripheral edge of a bacterial colony 
on an agar medium filling a Petri dish. The area around the edge is the theater of fast motion 
of swarmer cell flagella making it moves while the medium is gradually colonized.
 The last panel (D) shows a colony after several hours of swarming.  It forms a thin film  showing 
 the specific pattern: concentric rings called terraces.
Those rings are made by variations in the colony's thickness.
 
Figure \ref{figilsw} which
displays three snapshots of the colony expansion on agar during an 
active migration phase. The expansion process is clearly seen from 
the progression of the colony's edge. The central inoculum is 
located in the bottom right of each picture.

\:

For a precise and synthetic description of \emph{Proteus mirabilis} swarm phenomenon,
we refer to Rauprich \emph{et al.}  \cite{Rauetal1996} and  
Matsuyama \emph{et al.} \cite{MatTaItWaMat}. Among all aspects described in those references,
periodic alternation of swarm phases and of consolidation phases, without any motion,
seem to be the most typical one.
This alternation gives rise to terraces as a resulting effect.
The second one concerns the cell size distribution and was pointed out in 
Matsuyama \emph{et al.} \cite{MatTaItWaMat}. As for space distribution at a
given time, the swarming front is essentially constituted by long cells and the proportion 
of shorter cells grows while moving towards the colony interior. Concerning time 
distribution at a given place, the colonization  starts with long 
hyperflagellated cells and then the proportion 
of shorter cells grows.
The third typical fact, pointed out in  Lahaye \emph{et al.} \cite{LaAuKerDouSi,LaAuFleSi}
 is that the consolidation phase begins even though 
swarmer cells are still present in the colony front.

\: 

\:

\emph{Proteus mirabilis} swarm phenomenon has fascinated micro-biologists
and other scientists
for over a century, since Hauser \cite{Hau1885} first described it.
This fascination certainly comes from the fact that bacteria are often
considered as very primitive organisms and from the astonishing fact that 
such primitive living forms are able to generate such nice regular structures.
This phenomenon certainly played an important role in morphogenesis and fractal
theories (see Thom \cite{thom1}, Mandelbrot \cite{mand}, Falconer \cite{falc} 
and references in them) 
that highlighted in a
very resounding way that very simple principles may 
generate very complex shapes.
On the other hand, as explained by Kaufmann
\cite{kauf}, phenomena of the \emph{Proteus mirabilis} swarm kind
seem to belong to an intermediary step
between unicellular organisms and multicellular organisms. 
As a mater of fact, key elements to the understanding of evolution of life  
can certainly be found within the understanding of those phenomena.

Collective behaviour of bacteria colonies has been widely observed,
described and studied.
For instance, Hoeniger \cite{Hoe} delineated the \emph{Proteus mirabilis} 
swarm phenomenon and Budrene \& Berg \cite{BudBer} pointed out the ability 
that bacteria colonies have to draw complex patterns.
Shapiro \cite{Sha1985a,Sha1985d,Sha1987a} and more recently 
Rauprich \emph{et al.}  \cite{Rauetal1996} and  
Matsuyama \emph{et al.} \cite{MatTaItWaMat}  made a wide observation campaign
to precisely analyze the phenomenon.

This research effort pointed out that the swarm phenomenon has to be seen
as an emergent property leading to the consideration of bacteria as multicellular
organisms (see Shapiro \& Dworkin \cite{ShaDwo1997}, Shapiro \cite{Sha1988b,Sha} and
 Dworkin \cite{Dwo}).

\:

Research has mainly tackled the problems of origin and of regulation
of bacterial collective behaviour with a genetic or a bio-chemical slant
(see Shapiro \cite{Sha1985c,Sha1985b,Sha1993b,Sha1994a},
Belas \cite{Bel},
Manos \& Belas \cite{ManArtBe}, 
Schneider \emph{et al.} \cite{ScLoJoBe} and  Gu\'e \emph{et al.} \cite{Gueetal2001}) .

\:

\:

The goal of the present paper is to add a bio-physical slant concerning
the swarm phenomenon.
More precisely, we exhibit the main physical phenomena of the bacteria
to achieve their colonization. 
What we have chosen to do is to make assumptions about the physical principles
governing the bacteria behaviour and controlling the swarm
and to translate those assumptions into 
a mathematical model involving partial differential equations. Then,  
this model is implemented in Matlab$^\text{\textregistered}$7 language. Using the resulting
program, we make several simulations with several values of parameters involved
in the model. We finally observe that macroscopic colony behaviours 
may be recovered by the simulations
leading to the conclusion that our assumptions are not entirely unreasonable.
Beyond this conclusion, the scientific approach consisting in using modelling
and simulation in order to test assumptions, is another important contribution 
of this paper.
We plan to apply it, in forthcoming works, 
in order to refine \emph{Proteus mirabilis} swarm
explanation.

\:

Using modelling, mathematical analysis and simulation to describe or understand
biological systems is not unusual nowadays. For an introduction
to these topics, we refer to Murray \cite{mur}.
Those kinds of mathematical methods were applied in the context 
of collective behaviour of bacteria. For a description of this, we refer to 
Matsuschita \cite{mat} and 
Ben-Jacob \& Cohen \cite{BenjaCo}.

Esipov \& Shapiro \cite{EsiSha} set out a kinetic model describing the swarm of
\emph{Proteus mirabilis} at the colony scale recovering the main behaviours
of a real colony: time alternation of swarming phases and consolidation phases, 
time alternation of colony radius increases and stops and terraces. 
The key ingredient causing their
model to produce good results is the introduction
of an age structuring of the bacterial population as is suggested in
Gurtin \cite{Gur} and Gurtin \& Mac Camy \cite{GurMcCa1974}. This introduction 
allows one to take into consideration the age of 
the swarmer cells when looking at their behaviour.
This model was studied from mathematical and numerical analysis viewpoints in
Esipov \& Shapiro \cite{EsiSha}, Medvedev, Kapper \& Koppel \cite{MedKaKo},
Ayati \cite{Ayati2000,Ayati2005,Ayati2007,Ayati2008},
Ayati \& Dupont \cite{AyDu}, Fr\'enod \cite{fre2006},
Ayati  \& Klapper \cite{AyKla}
and Lauren\c{c}ot   \& Walker  \cite{LaurencotWalker2008, LaurencotWalker2009}.
In particular, Ayati \cite{Ayati2007} showed the importance
of using age structured model to reproduce the cell size distribution pointed out by  
Matsuyama \emph{et al.} \cite{MatTaItWaMat}.

The existing models describing colonization by bacteria generally contain 
an advection or diffusion
term allowing one to define the direction in which the colony expands 
conformal to reality.
This makes the model behave in a way which conforms to reality.
But the factors in front of these advection and diffusion terms which control 
the beginning, velocity, and ending of the swarm are mainly set from phenomenological
considerations. 

\: 

As for the research of physical principles explaining bacterial collective 
behaviour we mention Mendelson, Salhi \& Li \cite{MenSaLi}.

\:

\:

{\bf Acknowledgements - } We thank James Shapiro and  Bruce Ayati who pointed out 
lacks in the first version of the paper and suggested ways to remedy them.
We also thank Joanna Ropers for proofreading this paper. 

\section{Explanation of the principles governing the swarm}
\label{expl}
From a long \emph{in situ} observation campaign of \emph{Proteus mirabilis}
colonies (see
Gu\'e \emph{et al.} \cite{Gueetal2001}, 
Keirsse \emph{et al.} \cite{Keiretal},
Lahaye \emph{et al.} \cite{LaAuKerDouSi,LaAuFleSi}) and 
knowledge of physical phenomena arising within complex fluids, 
we are convinced that the main physical explanations of 
\emph{Proteus mirabilis} swarm come from the properties of the
extra cellular matrix which is the complex fluid smothering the bacteria of
the colony.
More than ten years ago, Rauprich \emph{et al.}  \cite{Rauetal1996} hypothesised
the key role of relative osmotic activities at the agar colony interface. More recently,
Berg \cite{berg} and Chen \emph{et al.} \cite{CheTuBe} further supported hypothesis of 
Lahaye \emph{et al.}  \cite{LaAuKerDouSi,LaAuFleSi}  claiming that the preservation 
of a critical water concentration allowing swarming of bacteria colonies was due to 
an osmotic agent.
  
In short,
two properties seem to play an important role: 
the thickness of the colony and the water concentration 
within the extra cellular matrix.

The colony thickness plays a role in controlling the growth rate of cells in
their vegetative form: once a given thickness is reached, 
the cell division ends, stopping vegetative population growth. 

The extra cellular matrix is compounded with water, polysaccharides and lipids. When
the water concentration of this matrix is low, polysaccharides and lipids 
auto-organize which results in semi-crystallinity and visco-elastic properties.
When it is high,  polysaccharides and lipid auto-organization is lost, 
due to efficient solvent competition with solute-solute interactions, 
decreasing thereby the visco-elastic properties of the extra cellular matrix.
\newline
It has been well established that the translocation ability of the colony
is a consequence of the rotation and oscillation the hundreds of  flagella of a large
number of swarmer cells. On the other hand, observations have put into light  that
consolidation phases begin despite  
numbers of swarmer cells are present and moving but without any effect on 
colony displacement.
\newline
Hence, we may infer that water concentration rules the swarm phenomenon
in the following way. When it is low, swarmer cells interacting with visco-elastic
properties of the extra cellular matrix produce colony expansion. 
When the matrix water concentration is high, the visco-elastic properties are lost. 
Then swarmer cell agitation, when it occurs, does not produce a swarm.

The differentiation process may be assumed to be a kind of dysfunction in
the cell division process: for a small proportion of vegetative bacteria,
the cell division is not brought to completion. It is also consistent
to suppose that a cell having undergone this 
differentiation goes on elongating, as long as water concentration is low,
until it reaches a maximum length depending on water concentration.
\newline
Once this maximum length is reached, the cell oscillation and rotation begin.
At this stage, the swarmer cells  become almost metabolically dormant
(see for instance Stickler \cite{stickler}, Mc Coy \emph{et al.} \cite{MCoLFaGu} and  Falkinham and Hoffman \cite{FalHo}) and then
do not consume extra fuel from the extra cellular matrix 
during their migration.
Then we may suppose that their life duration, in this stage,
is finite and linked to their length. After their life duration, 
the metabolism-frozen swarmer cells
de-differentiate to produce new vegetative cells. 

\:

In a given place of the colony, the  evolution of water concentration 
of the extra cellular matrix is influenced by the density of bacteria, which
need water to achieve metabolism, and by the displacement of swarmer cells which
bring with them a part of their surrounding medium. 
The third factor influencing matrix water concentration is
the water concentration in the agar medium with
which it exchanges water and the corresponding transfer properties 
at the agar - colony interface.

\:

\:

\section{Translation of the principles into a model} 
\label{trans}
In this section, we translate the principles
explained in the previous section into a model,
in the simplest manner possible,
and even in an oversimplified manner.

\:

The model set out in  Esipov \& Shapiro \cite{EsiSha} coupled
two equations. The first of them was an ordinary differential
equation describing the evolution of the vegetative cell biomass density.
The swarmer cells were modelled by an age-structured density which was
the solution to a partial differential equation describing aging, but also 
de-differentiation and swarming, both being age dependent processes.

The model we have set out in the present paper, involves three equations.
The first concerns vegetative cells. 
To be able to take into account that the differentiated 
cells undergo an elongation phase, followed by a swarming phase 
with no metabolism,
we consider an equation modelling the evolution of the elongating cells 
and another one modelling the evolution of swarmer cells having stopped
their metabolism, both of them involving age-structured functions.

\:

\:

We begin by considering the evolution of the vegetative cells. For this,
we introduce the function
$\qc(\tc,\xc)$ which is the vegetative cell biomass density at time 
$t\in[0,T)$ and in a given position $\xc\in \ritd$.
It is supposed to be the solution to the following ordinary differential
equation:
\begin{gather}
\label{S1}
\ds \fracp {\qc}{\tc} 
= \frac{1-\xic}{\tac}\: \qc \; \chc\Big(\frac{\bioz-\biom}{\bioz} \Big)
\; \chc\Big(\frac{\qc-\qcz}{\qcz} \Big)
+ \int_0^{+\infty} \rhc(\tc,\age,\kac \age,\xc) \: e^{\age/\tac} \, d\age.
\end{gather}
The time evolution of $\qc$ is the result of cell division, described
by the first term on the right-hand side of the equation (\ref{S1}), and of
the de-differentiation products quantified by the second term.
In the first term, $\tac$  stands for the growth rate of the cell division process,
 $\xic$, $0\leq \xic\leq 1$, is the proportion of cells undergoing 
differentiation, and then
$1-\xic$ is the proportion of cells
 bringing their division to completion.
The function $\chc$ has the following definition
$\chc(E)=0$ if $E<0$ and $=1$ otherwise.
Then the factor $\chc((\bioz-\biom)/\bioz)$, where $\biom= \biom(\tc,\xc)$ is 
swarm colony thickness at $\tc$ in $\xc$ , takes into account the above evoked control of
colony thickness on the cell division process in the following way: when
$\biom\geq \bioz$, the term is set to $0$. 
The factor $\chc((\qc-\qcz)/\qcz)$ precludes the growth of $\qc$
by cell division when $\qc$ is too low, \emph{i.e.} when $\qc<\qcz$.  
\newline
In the second term, 
$\rhc(\tc,\age,\bge,\xc)$, at time $\tc\in[0,T)$, stands for the 
density of swarmer cells in position $\xc\in \ritd$, having been freezing
their metabolism at an age $\age\in[0,+\infty)$, for a period of time 
$\bge$. Since the life duration of the the swarmer cells in this state is
supposed to be linked to their length and since the length of a given swarmer
cell is linked to the duration during which it elongated,
we set that the life duration of a swarmer cell
having stopped metabolism at an age $\age$ is $\kac\age$. 
The growth rate of the biomass of the elongating cells 
is the same as the one of the vegetative cells.
Hence the biomass density in position $\xc$ of swarmer cells  
having been stopping 
their metabolism at an age $\age$, for a period of time 
$\bge$ is $\rhc(\tc,\age,\bge,\xc) \: e^{\age/\tac}$ and the integral
in the second term on the right-hand side
 of (\ref{S1}) is nothing but
the biomass density of cells that de-differentiate at time $\tc$.

The biomass density of vegetative cells at time $\tc=0$ is given by 
the following initial data: 
\begin{gather}
\label{S2}
\qc(0,\xc) = \qcinit(\xc).
\end{gather}

\:

\:

The description of the evolution of elongating cell population
involves the following age-structured function 
$\zec(\tc,\age,\xc)$, which stands, at time $\tc\in[0,T)$, 
for the density of elongating cells of age $\age\in[0,+\infty)$ in $\xc\in \ritd$.
The evolution of this density is due to aging and to the transformation of 
elongating cells into swarmer cells. We assume that this transformation 
in a given position $\xc$
takes place according to a transformation rate $\muc(\age,\hhc)$ depending
on age $\age$ of the considered bacteria and on water concentration 
$\hhc=\hhc(\tc,\xc)$ of the matrix. The dependence of $\muc$ with respect to 
$\hhc$ takes into  account that the lower $\hhc$ is, the larger the existence
duration expectancy of elongating cells is.

The equation modelling the evolution of $\zec$ is:
\begin{gather}
\label{S3}
\fracp{\zec}{\tc} + \fracp{\zec}{\age}
= - \muc(\age,\hhc) \zec.
\end{gather}

Since the proportion $\xic$ of vegetative cells that differentiate
generates elongating cells with age $0$, equation (\ref{S1}) and  (\ref{S3})
are linked by the following transfer condition,
\begin{gather}
\label{S4}
\zec(\tc,0,\xc) = \frac{\xic}{\tac}\: \qc(\tc,\xc) \; \chc\Big(\frac{\bioz-\biom}{\bioz} \Big)
\; \chc\Big(\frac{\qc-\qcz}{\qcz} \Big).
\end{gather}

The elongating cell density at the initial time is given by:
\begin{gather}
\label{S5}
\zec(0,\age,\xc) = \zec_0(\age,\xc),
\end{gather}

\:

\:

Now we turn to the modelling of the evolution of the swarmer cells
which have stopped their metabolism. 
\newline
We consider that the velocity of swarmer cells is the consequence of the
interaction  of their agitation with the visco-elastic properties of the
matrix when water concentration is low. 
Then we introduce a function
$c(\hhc)$ which is non-zero for the small value of $\hhc$ and with value $0$
for the large value of  $\hhc$. Then, as the global motion of swarmer cells, 
when it exists, is the consequence of cell agitation it is reasonable to consider
that their velocity is in proportion with the colony thickness gradient
$\nabla \biom$, the proportion coefficient being $c(\hhc)$. In other words,
we consider that the swarmer cell velocity in a given position $\xc$ is 
$c(\hhc(\tc,\xc)) \nabla \biom(\tc,\xc)$.
\newline
The definition of 
$\rhc(\tc,\age,\bge,\xc)$ defined for $t\in[0,T)$, $\age\in[0,+\infty)$,
$\bge\in[0,\kac\age)$ and $\xc\in \ritd$ was already given.
In view of this definition, if $S$ is a regular set of $\ritd$
and $\bge_1 < \bge_2\leq \kac \age$ are two numbers, 
\begin{gather}
\ds \int_{S\times[\bge_1,\bge_2]} \rhc(\tc,\age,\bge,\xc)\, d\bge d\xc,
\end{gather}
is, at time $\tc$, the density of bacteria having stopped their metabolism at an age $\age$,
for a time larger than $\bge_1$ and smaller than $\bge_2$ and situated in $S$.
The time evolution of this quantity is linked to the aging fluxes in $\bge_1$ and 
$\bge_2$ and to the swarming flux on the boundary $\partial S$ of $S$ which has the 
following expression
\begin{gather}
\int_{\partial S \times[b_1,b_2]} \big(c(\hhc) \nabla \biom \: \rhc\big) \cdot n \, d\bge dl,
\end{gather}
where $n$ stands for vector with norm $1$, orthogonal to $S$ and pointing outwards
and where $dl$ is the Lebesgue measure on $\partial S$. 
Hence we deduce the following equation:
\begin{gather}
\label{S66}
\frac{\ds d\Big(\int_{S\times[\bge_1,\bge_2]} \rhc\, d\bge d\xc \Big)}{dt}
= \int_{S}\rhc(\tc,\age,\bge_1,\xc)\, d\xc 
- \int_{S}\rhc(\tc,\age,\bge_2,\xc)\, d\xc
- \int_{\partial S \times[b_1,b_2]} \big(c(\hhc) \nabla \biom \: \rhc\big) \cdot n \, d\bge dl.
\end{gather}
Using the divergence operator $\nabla\cdot$ to transform flux integral into
integral over $S$ and noticing that 
$\nabla \cdot (c(\hhc) \nabla \biom \: \rhc)
= c(\hhc) \nabla \biom \cdot \nabla\rhc + \rhc \,\nabla\cdot (c(\hhc) \nabla \biom)$,
we finally obtain the partial differential equation $\rhc$ satisfies:
\begin{gather}
\label{S6}
\fracp{\rhc}{\tc} + \fracp{\rhc}{\bge}
= -c(\hhc) \nabla \biom \cdot \nabla\rhc - \rhc \,\nabla\cdot(c(\hhc)\nabla \biom).
\end{gather}
Since the elongating cells that stop metabolism at an age $\age$
become swarmer cells, we have the following transfer condition
\begin{gather}
\label{S7}
\rhc(\tc,\age,0,\xc)= \muc(\age,\hhc)\zec,
\end{gather}
and the following initial condition
\begin{gather}
\label{S8}
\rhc(\tc,0,\bge,\xc)=\rhc_0(\tc,\bge,\xc),
\end{gather}
to provide (\ref{S6}) with.

\:

\:

The thickness $\biom$ of the colony in a given position is in direct
proportion to the biomass density value at this point.
In other words we set
\begin{gather}
\label{defbiom}
\biom = \qc + \biomm + \biomn,
\end{gather}
where 
\begin{gather}
\label{defbiomm}
\biomm (\tc,\xc) = \int_0^{+\infty} \zec(\tc,\age,\xc)\; e^{\age/\tac} \, d\age,
\end{gather}
is the biomass density of elongating cells. Notice that  $\qc +\biomm$ is the
biomass density of metabolic cells; and $\biomn,$ defined as 
\begin{gather}
\label{defbiomn}
\biomn (\tc,\xc) 
= \int_0^{+\infty}\int_0^{\kac\age}\rhc(\tc,\age,\bge,\xc)\; e^{\age/\tac} \, d\age d\bge,
\end{gather}
is the density of metabolically dormant cells.

\:

\:

We now turn to the description of water concentration evolution.
We first consider the water concentration of agar gel near
the contact surface. Let
$\ggc(\tc,\xc)$, defined for $t\in[0,T)$, $\xc\in \ritd$, be 
the value of this water concentration at time $\tc$ 
and in position $\xc$. We consider that it evolves due to water transfers
from and to the extra cellular matrix and that those transfers
are proportional to the water concentration difference between 
the two media and to a colony thickness function. 
We name $\gamt$ the proportion factor.
There are also water exchanges between the agar layers close to the
contact surface and the deeper layers. We model those exchanges
by a simple ability to relax to the concentration
value being worth 1 with
a velocity $\gamd$. In other words, we consider that $\ggc$ evolves
according to the following equation:
\begin{gather}
\label{S9}
\fracp{\ggc}{t}= - \gamt\trnc(\biom) (\ggc-\hhc) + \gamd (1-\ggc),
\end{gather}
where $\trnc(\biom) = \biom$ if $\biom\leq1$ and $\trnc(\biom)=1$ otherwise.
Equation (\ref{S9}) is equipped with the following initial data
\begin{gather}
\label{S10}
\ggc(0,\xc)=\ggc_0(\xc),
\end{gather}
giving water concentration of agar at initial time.

\:

\:

Since it is easier to model water quantity transfer as water concentration
transfer, we introduce $\hc(\tc,\xc)$, defined for  $t\in[0,T)$ and $\xc\in \ritd$
which is the water density contained within the matrix in position $\xc$ at time 
$\tc$. This density is clearly proportional to the thickness of the colony times
the water concentration, with a proportion factor $\etac$ that is linked to the 
part of the swarm colony the matrix constitutes, \emph{i.e.}
$\hc = \etac\biom\hhc$.
\newline
If $S$ is a regular subset of $\ritd$ with boundary $\partial S$
equipped with the vector field $n$ with norm $1$ orthogonal to $S$ and pointing outwards
and with the Lebesgue measure $dl$, 
the quantity
\begin{gather}
\int_{S} \hc\, d\xc,
\end{gather}
is the water quantity contained in the swarm colony over $S$.
This quantity evolves according to
\begin{multline}
\label{S13}
\frac{\ds d\Big(\int_{S} \hc\, d\xc \Big)}{d\tc}
= - \int_{S} \Bigg(
\alc - \bigg(1-\chc\Big(\frac{\bioz-\biom}{\bioz}\Big)\,
               \chc\Big(\frac{\qc-\qcz}{\qcz} \Big)  \bigg) \alc'
             \Bigg) \qc \, d\xc 
- \int_{S} \alc \biomm \, d\xc 
\\
+ \int_{S} \gamt \trnc(\biom) (\ggc-\hhc) \, d\xc
-\int_{\partial S}\big( c(\hhc) \nabla\biom \, \etac \biomn \hhc\big) \cdot n \; dl.
\end{multline}
The second term of the right-hand side
 of this equation quantifies the amount
of water the elongating cells use for metabolism. This quantity is proportional, with
a factor $\alc$, to the biomass density  of elongating cells. 
The first term quantifies the amount of water used for metabolism by vegetative cells. 
This quantity is proportional to the biomass density of vegetative cells.
The proportion factor is $\alc$ when vegetative cells undergo cell division 
and is $\alc -\alc'$ when they do not undergo cell division, \emph{i.e.}
when $\chc((\bioz-\biom)/\bioz)\chc((\qc-\qcz)/\qcz)=0.$
Since the swarmer cells are metabolically dormant, we consider they do not 
intervene in water consumption.
The third term models the water exchanges with agar.
The last term gives the flux of water carried by the swarmer cells. In it,
$c(\hhc) \nabla \biom $ is the velocity of the swarmers and then of the medium 
they carry with them; $\etac \biomn \hhc$ is the water quantity carried by the
swarmer cells, it is proportional to the swarmer biomass and to the water 
concentration in the matrix.

Using the divergence operator, we deduce from (\ref{S13}) that $\hc$ is solution to
\begin{multline}
\label{S14}
\fracp{\hc}{\tc} = - \Bigg(
\alc - \bigg(1-\chc\Big(\frac{\bioz-\biom}{\bioz}\Big)\,
               \chc\Big(\frac{\qc-\qcz}{\qcz} \Big)  \bigg) \alc'
                     \Bigg) \qc 
- \alc \biomm 
\\
+ \gamt  \trnc(\biom) (\ggc-\hhc) 
- \etac c(\hhc) \nabla\biom \cdot \nabla( \biomn \hhc) 
- \etac \biomn \hhc\,  \nabla\cdot(c(\hhc)\nabla\biom),
\end{multline}
endowed with the following initial data  
\begin{gather}
\label{S1414}
\hc(0,\xc)= \hc_0(\xc),
\end{gather}
describing the water density at time $\tc=0$.

From the link between $\hhc$ and $\hc$, we deduce 
\begin{gather}
\label{S15}
\fracp{\hc}{\tc} = \etac \Big(\hhc \fracp{\biom}{\tc} 
+ \biom \fracp{\hhc}{\tc}\Big),
\end{gather}
and then, from equation (\ref{S14}), an equation for $\hhc$
\begin{gather}
\nonumber
\fracp{\hhc}{\tc} = - \frac{\hhc}{\biom}\fracp{\biom}{\tc}  
- \frac{1}{\etac} \Bigg(
\alc - \bigg(1-\chc\Big(\frac{\bioz-\biom}{\bioz}\Big)\,
               \chc\Big(\frac{\qc-\qcz}{\qcz} \Big)  \bigg) \alc'
                     \Bigg)\frac{\qc}{\biom} ~~~~~~~~~~~~~~~~~~~~~~~~~~~~~~~
\\
\label{S16}
- \frac{\alc}{\etac} \frac{\biomm}{\biom}
+ \frac{\gamt \trnc(\biom)}{\etac\biom} (\ggc-\hhc) 
- \frac{c(\hhc)}{\biom} \nabla\biom \cdot \nabla(\biomn \hhc) 
- \frac{\biomn}{\biom} \hhc\, \nabla\cdot(c(\hhc)\nabla\biom),
\\
\label{S1616}
\hhc(0,\xc)= \hhc_0(\xc) = \frac{\hc_0(\xc)}{\etac\biom}.
\end{gather}

\section{Implementation}
\label{impl}
From the model just set out, we write a Biomass Preserving Finite Volume-like 
numerical scheme (see LeVeque \cite{LeV1992,LeV2002}) and then we implement it.

\:

\:

\begin{figure}
\setlength{\unitlength}{1mm}
\centering
\begin{picture}(95,20)(0,10)
\put(0,15){\line(1,0){100} }
\put(99,15.0){\makebox(0,0){$>$}}
\put(99.3,15.0){\makebox(0,0){$>$}}
\put(7,15){\makebox(0,0){$\bullet$}}
\put(7,12){\makebox(0,0){$\age_{1/2}$}}
\put(7,20){\makebox(0,0){$0$}}
\put(14,15){\makebox(0,0){$|$}}
\put(14,10){\makebox(0,0){$\age_{1}$}}
\put(21,15){\makebox(0,0){$\bullet$}}
\put(21,12){\makebox(0,0){$\age_{3/2}$}}
\put(21,20){\makebox(0,0){$\leftarrow \agest\rightarrow$}}
\put(28,15){\makebox(0,0){$|$}}
\put(28,10){\makebox(0,0){$\age_{2}$}}
\put(35,15){\makebox(0,0){$///$}}
\put(43,15){\makebox(0,0){$\bullet$}}
\put(43,12){\makebox(0,0){$\age_{k-1/2}$}}
\put(50,15){\makebox(0,0){$|$}}
\put(50,10){\makebox(0,0){$\age_{k}$}}
\put(57,15){\makebox(0,0){$\bullet$}}
\put(57,12){\makebox(0,0){$\age_{k+1/2}$}}
\put(67,15){\makebox(0,0){$///$}}
\end{picture}
\caption{Elongating cell age axis discretization}\label{figadisc}
\end{figure}
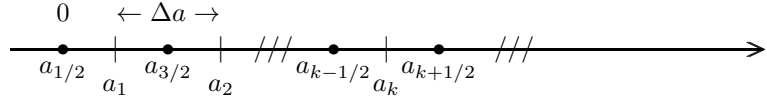

\begin{figure}
\setlength{\unitlength}{1.2mm}
\centering
\begin{picture}(95,20)(0,10)
\put(0,15){\line(1,0){95} }
\put(94,15.0){\makebox(0,0){$>$}}
\put(94.3,15.0){\makebox(0,0){$>$}}
\put(7,15){\makebox(0,0){$\bullet$}}
\put(7,12){\makebox(0,0){$\xc_{1/2}$}}
\put(7,20){\makebox(0,0){$0$}}
\put(14,15){\makebox(0,0){$|$}}
\put(14,10){\makebox(0,0){$\xc_{1}$}}
\put(21,15){\makebox(0,0){$\bullet$}}
\put(21,12){\makebox(0,0){$\xc_{3/2}$}}
\put(21,20){\makebox(0,0){$\leftarrow \, \posst \,\rightarrow$}}
\put(28,15){\makebox(0,0){$|$}}
\put(28,10){\makebox(0,0){$\xc_{2}$}}
\put(35,15){\makebox(0,0){$///$}}
\put(43,15){\makebox(0,0){$\bullet$}}
\put(43,12){\makebox(0,0){$\xc_{i-1/2}$}}
\put(50,15){\makebox(0,0){$|$}}
\put(50,10){\makebox(0,0){$\xc_{i}$}}
\put(57,15){\makebox(0,0){$\bullet$}}
\put(57,12){\makebox(0,0){$\xc_{i+1/2}$}}
\put(64,15){\makebox(0,0){$///$}}
\put(72,15){\makebox(0,0){$\bullet$}}
\put(72,12){\makebox(0,0){$\xc_{I-1/2}$}}
\put(79,15){\makebox(0,0){$|$}}
\put(79,10){\makebox(0,0){$\xc_{I}$}}
\put(86,15){\makebox(0,0){$\bullet$}}
\put(86,12){\makebox(0,0){$\xc_{I+1/2}$}}
\put(86,20){\makebox(0,0){$\posmax$}}
\end{picture}
\caption{Position axis discretization}\label{figxdisc}
\end{figure}
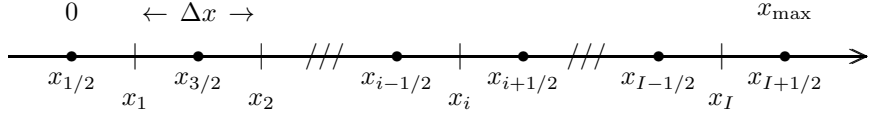

In this section, to simplify matters, we consider the swarm of a 
\emph{Proteus mirabilis} population whose distribution remains independent 
of the second position variable and an agar medium which is also independent 
of the second position variable.
This means that the position space is considered nothing but
a segment, say $[0,\posmax]$, of $\rit$.
We introduce a time step $\timest$, and age steps $\agest$
and $\bgest$ with the following link $\bgest=\agest=\timest$.
We also define a position step $\posst$ such that,
$I \posst = \posmax$, for an integer $I$.
We then define, for $n\in \nit$,  $\tc_n = n \timest$,
for $k \in \nit^*$, $\age_k= (k-1/2) \agest$ (see Figure \ref{figadisc}),
for $p\in \nit^*$, $\bge_ p= (p-1/2) \bgest$, and, 
for $i \in \{1,\dots, I\}$, $\xc_i = (i-1/2) \posst$
(see Figure \ref{figxdisc}).

\:

For every $\tc_n$, the function $\qc(\tc_n,.)$ is approximated 
by a function with a constant value, defined as $\qc_i^n$, on every interval
$(\xc_{i-1/2},\xc_{i+1/2})$, where $\xc_{i-1/2}=\xc_i- \posst/2= (i-1)\posst$ and
$\xc_{i+1/2}=\xc_i+ \posst/2= i\posst$ for $i\in \{1,\dots, I\}$.
Consistently, the approximation of initial data (\ref{S2}) is:
\begin{gather}
\label{A1}
\qc_i^0 = \frac{1}{\posst} \int_{x_{i-1/2}}^{x_{i+1/2}}\qcinit(\xc) \, d\xc.
\end{gather}
Equation (\ref{S1}) is approximated by 
\begin{gather}
\label{A2}
\qc_i^{n+1} = \qc_i^n
+ \timest \bigg(
  \frac{1-\xic}{\tac} \qc_i^n\;\chc\Big(\frac{\bioz-\biom_i^n}{\bioz} \Big)
  \, \chc\Big(\frac{\qc_i^n-\qcz}{\qcz} \Big)
+ \rrc_i^n
          \bigg),
\end{gather}
where $\biom_i^n$ and $\rrc_i^n$ are approximations of the mean values 
of $\biom(\tc_n,.)$ and of 
$\int_0^{+\infty} \rhc(\tc_n,\age,\kac \age,.) \: e^{\age/\tac} \, d\age$ over
$(\xc_{i-1/2},\xc_{i+1/2})$, supposed to be available and defined in the sequel.

\:

\begin{figure}
\centering 
\includegraphics[height=5.5cm]{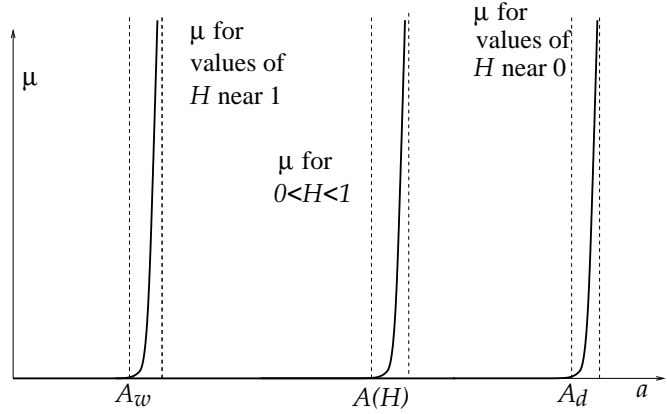}
\caption{Shapes of $\muc$ for matrix with high, medium or low water concentration.}
\label{figAofH}
\end{figure}

The function $\zec(\tc_n,.,.)$  is then approached by a function
$\tilde \zec(\tc_n,.,.)$, which is constant, with worth $\zec^{n}_{k,i}$, on any 
rectangle 
${\cal R}_{k,i} = (\age_{k-1/2},\age_{k+1/2})\times (\xc_{i-1/2},\xc_{i+1/2})$,
where $\age_{k-1/2}=\age_k-\agest/2=(k-1)\agest$ and 
$\age_{k+1/2}=\age_k+\agest/2=k\agest$, such that 
\begin{multline}
\label{A21}
\int_{{\cal R}_{k,i} }
\tilde \zec(\tc_n,\age,\xc) \, e^{\age/\tac} d\age d\xc
=
\posst \,\big(e^{ k\agest/\tac}-e^{(k-1)\agest/\tac}\big)\tac\; \zec^{n}_{k,i} 
\\
\text{is supposed to be close to}
\int_{{\cal R}_{k,i} }
\zec(\tc_n,\age,\xc) \, e^{\age/\tac}\; d\age d\xc,
\end{multline}
for $k\in \nit^*$ and $i\in \{1,\dots, I\}$. Then, the approximation of
initial data (\ref{S5}) reads
\begin{gather}
\label{A3}
\zec^{0}_{k,i} = \frac{1}{\posst \,(e^{k\agest/\tac}-e^{(k-1)\agest/\tac})\tac}
\int_{{\cal R}_{k,i} }
\zec_0(\age,\xc) \, e^{\age/\tac} \; d\age d\xc.
\end{gather}
At each time step $\tc_n$, $\zec^{n}_{0,i}$ for $i\in \{1,\dots, I\}$
are used for the biomass transfer from vegetative cells. Hence we translate
condition (\ref{S4}) into:
\begin{gather}
\label{A4}
\zec^{n}_{0,i}= \frac{\timest}{~(e^{\agest/\tac}-1)\tac~}\;
  \frac{\xic}{\tac}\; \qc_i^n \; \chc \Big(\frac{\bioz-\biom_i^n}{\bioz} \Big)
   \, \chc\Big(\frac{\qc_i^n-\qcz}{\qcz} \Big).
\end{gather}
We introduce a non increasing function $A(H)$, ranging in the interval
$[A_w, A_d]$ with $A(0)=A_d$ and $A(1)=A_w$, and 
before approximating (\ref{S3}) we assume that for any value of $H$,
function $\muc(.,H)$ is 0 when  $a\leq A(H)$ and increases very quickly after this
value, as illustrated in Figure \ref{figAofH}. 
In other words, we assume that at a given value $H$ of the water concentration,
the elongating cells remain in there state
 until an age $A(H)$ and that
they become swarmer cells very soon after this age.
With this assumption and since $\agest=\timest$, equation (\ref{S3}) may be approximated by:
\begin{align}
\label{A5}
\zec^{n+1}_{k,i} &= \zec^{n}_{k-1,i}, 
\text{ for every } i \in \{1,\dots, I\} 
\text{ and every } k \in \nit^* \text{ such that } (k-1/2)\agest \leq A(H^n_i), \\
\label{A6}
\zec^{n+1}_{k,i} &= 0,
\text{ for every } i \in \{1,\dots, I\} 
\text{ and every } k \in \nit^* \text{ such that } (k-1/2)\agest > A(H^n_i),
\end{align}
where $H^n_i$ is an approximation of $H$ near $\xc_i$ at time $\tc_n$ to be 
defined hereafter.

\:

We now turn to the swarmer cells. We approximate their density $\rhc(\tc_n,.,.)$
by a function $\tilde\rhc$ which is constant on any parallelepiped
${\cal P}_{k,p,i}=(\age_{k-1/2},\age_{k+1/2})\times(\bge_{p-1/2},\age_{p+1/2})
\times (\xc_{i-1/2},\xc_{i+1/2})$. The value $\rhc^{n}_{k,p,i}$ 
is defined as being such that $\bgest \, \rhc^{n}_{k,p,i}$ 
is the value of  $\tilde\rhc$ on this parallelepiped.
Consistent with this definition, initial data (\ref{S8}) is approximated by
\begin{gather}
\label{A7}
\rhc^{0}_{k,p,i} = \frac{1}{~\posst (e^{ k\agest/\tac}-e^{(k-1)\agest/\tac})\tac~}
\int_{{\cal P}_{k,p,i} }
\rhc_0(\age,\bge,\xc) \, e^{\age/\tac} \; d\age d\bge d\xc.
\end{gather}
Transfer condition (\ref{S7}) is approached by:
\begin{gather}
\label{starstar}
\rhc^{n+1}_{k,0,i} = \sum_{(k-1/2)\agest > A(H^n_i)} \zec^{n}_{k,i},
\end{gather}
for every $i \in \{1,\dots, I\}$, this last relation balancing
(\ref{A6}).
\newline
To be able to obtain a biomass preserving scheme, we chose to
discretize the form (\ref{S66}) of the equation giving the evolution of $\rhc$.
We make a splitting
 that consists in sequentially processing aging,
the action of which is quantified by the two first terms on the right-hand side
of equation (\ref{S66}), and, translocation modelled by the last term.
\newline
The first step of the splitting that takes into account aging, consists
simply in introducing a sequence $(\rhc^{n+1/2}_{k,p,i})$ for $k \in \nit^*$,
$p\in \nit^*$ and $i \in \{1,\dots, I\}$ defined by
\begin{gather}
\label{A8}
\rhc^{n+1/2}_{k,p,i} = \rhc^{n}_{k,p-1,i}.
\end{gather}
For the second step of the splitting
 we need an approximation
of the flux velocity $c(\hhc) \nabla \biom$
at all points $\xc_{i-1/2} = (i-1) \posst$
which are defined for $i \in \{1,\dots, I+1\}$, and which are the interfaces
between the discretization interval $(\xc_{i-1/2},\xc_{i+1/2})$ of the 
position space. Calling $V^n_{i-1/2}$ the value of this approximation, we set
$i \in \{2,\dots, I\}$
\begin{gather}
\label{A9}
V^n_{i-1/2} = c\Bigg(\frac{2}{\frac{1}{\hhc^{n}_{i-1}}+\frac{1}{\hhc^{n}_{i}}}\Bigg)
\, \frac{1}{\posst} (\biom^{n}_{i}-\biom^{n}_{i-1}),
\end{gather}
that allows one to
minimize retrograde polluting propagation of the flux in 
regions where its velocity is zero. We set $V^n_{1/2}=0$ and 
$V^n_{I+1/2}=V^n_{I-1/2}$ which makes the swarmer cells go out of the domain boundary once they
have reached it.
\newline
From those flux velocities, we define the swarmer fluxes $F^{n}_{k,p,i-1/2}$ in all interfaces
$\xc_{i-1/2}$. In an interface where the velocity is non negative, the 
value of the swarmer density just before the interface is taken into account.
In an interface where the velocity is non positive, the 
value of the swarmer density just after the interface is taken into account.
In other words we set:
\begin{equation}\label{A12}\begin{aligned}
F^{n}_{k,p,i-1/2} &= V^n_{i-1/2} \, \rhc^{n+1/2}_{k,p,i} \text{ ~ if ~ } V^n_{i-1/2} \leq 0,
\\
F^{n}_{k,p,i-1/2} &= V^n_{i-1/2} \,\rhc^{n+1/2}_{k,p,i-1} \text{ ~ if ~ } V^n_{i-1/2} > 0.
\end{aligned}\end{equation}
Once the fluxes are computed, not forgetting that $\rhc$ is defined
for values of $\bge$ smaller than $\kac \age$ and that swarmer cells at an age
$\bge=\kac \age$ de-differentiate to produce vegetative cells,
the values of the swarmer density approximation for the next time step is given,
for every $i \in \{1,\dots, I\}$, by
\begin{align}
\rhc^{n+1}_{k,p,i} = &\rhc^{n+1/2}_{k,p,i} + \frac{\timest}{\posst}
(F^{n}_{k,p,i-1/2} - F^{n}_{k,p,i+1/2} ),
\nonumber
\\
\label{A13}
& ~~~~~~~~~~~~~~~~~ 
\text{for every } p \in \nit^* \text{ and } k \in \nit^* 
\text{ such that } (p-1/2)\bgest\leq \kac (k-1/2)\agest,
\\
\label{A13.1}
\rhc^{n+1}_{k,p,i} = &0,
\text{for every } p \in \nit^* \text{ and } k \in \nit^* 
\text{ such that } (p-1/2)\bgest > \kac (k-1/2)\agest,
\end{align}
and the approximation $\rrc_i^{n+1}$ of 
$\int_0^{+\infty} \rhc(\tc_{n+1},\age,\kac \age,.) \: e^{\age/\tac} \, d\age$
available for the next time step, by 
\begin{gather}
\label{A13.2}
\rrc_i^{n+1} = \sum_{(p-1/2)\bgest > \kac (k-1/2)\agest}
\bigg(\rhc^{n+1/2}_{k,p,i} + \frac{\timest}{\posst}
(F^{n}_{k,p,i-1/2} - F^{n}_{k,p,i+1/2} )\bigg)
\big(e^{k\agest/\tac} -e^{(k-1)\agest/\tac}\big)\tac.
\end{gather}

\:

In every formula above the approximation $(\biom^n_i)$ of the thickness
$\biom$ is given, for  $i\in \{1,\dots, I\}$, by
\begin{gather}
\label{A15}
\biom^n_i =  \qc^n_i + \biomm^n_i  + \biomn^n_i,
\\
\label{A151}
\biomm^n_i = \sum_{k \in \nit^*} 
\big(e^{k\agest/\tac}-e^{(k-1)\agest/\tac}\big)\tac\;
\zec^{n}_{k,i},
\\
\label{A152}
\biomn^n_i = \sum_{(p-1/2)\bgest\leq\kac(k-1/2)\agest} 
\big(e^{k\agest/\tac }-e^{(k-1)\agest/\tac}\big)\tac\; \rhc^n_{k,p,i},
\end{gather}
in accordance with (\ref{defbiom}), (\ref{defbiomm}) and (\ref{defbiomn}).
 
Function $\ggc(t_n,\cdot)$ is approximated by a function which is 
constant on each interval $(\xc_{i-1/2},\xc_{i+1/2})$.
Equation (\ref{S9}) and initial data (\ref{S10}) are approached by
\begin{gather}
\label{A17000}
\ggc^{n+1}_{i}= \ggc^{n}_{i}
- \timest\,\gamt\trnc(\biom^n_i) (\ggc^{n}_{i}-\hhc^{n}_{i}) 
+ \timest\,\gamd (1-\ggc^{n}_{i}),
\\ \ds
\label{A17AAA}
\ggc^{0}_{i} = \frac{1}{\posst} \int_{x_{i-1/2}}^{x_{i+1/2}} \ggc_0(\xc) \; dx.
\end{gather}

\:

Finally, we approximate the evolution of the water quantity by 
\begin{multline}
\label{A1600}
\hc^{n+1}_{i} = \hc^{n}_{i} 
- \timest\, \Bigg(
\alc - \bigg(1-\chc\Big(\frac{\bioz-\biom^{n}_{i}}{\bioz}\Big)\,
               \chc\Big(\frac{\qc^{n}_{i}-\qcz}{\qcz} \Big)  \bigg) \alc'
             \Bigg) \qc^{n}_{i}
- \timest\, \alc \biomm^{n}_{i}
\\
+ \timest\,\gamt \trnc(\biom^{n}_{i})(\ggc^{n}_{i}-\hhc^{n}_{i})
+ \frac{\timest}{\posst}\,(\tilde F^{n}_{i-1/2}-\tilde F^{n}_{i+1/2})
\end{multline}
where, $V^n_{i-1/2}$ being defined by (\ref{A9}),
\begin{equation}\label{A17}\begin{aligned}
\tilde F^{n}_{i-1/2} &= \etac\: V^n_{i-1/2} \, \biomn^{n}_{i} \hhc^{n}_{i}
\text{ ~ if ~ } V^n_{i-1/2} \leq 0,
\\
\tilde F^{n}_{i-1/2} &= \etac\: V^n_{i-1/2} \,\biomn^{n}_{i-1} \hhc^{n}_{i-1}
\text{ ~ if ~ } V^n_{i-1/2} > 0,
\end{aligned}\end{equation}
and where $\hc^{n}_{i}$ is the value of the constant by interval approximation
of $\hc(\tc_n,\cdot)$ on $(\xc_{i-1/2},\xc_{i+1/2})$.
\newline
We then get
 $\hhc^{n+1}_{i}$ by 
\begin{gather}
\label{A17BBB}
\hhc^{n+1}_{i} = \frac{\hc^{n+1}_{i}}{\etac \,\biom^{n+1}_{i}}.
\end{gather}
The scheme just built was implemented in Matlab$^\text{\textregistered}$7 language.

\section{Simulations}
\label{simul}
In this section we present several simulations for various values 
of the given parameters
Those simulations lead to the conclusion that 
characteristics of the \emph{Proteus mirabilis} colony can be recovered
by the model just built.

\subsection{Simulation 1}
\begin{figure}
\centering 
\includegraphics[width=16cm]{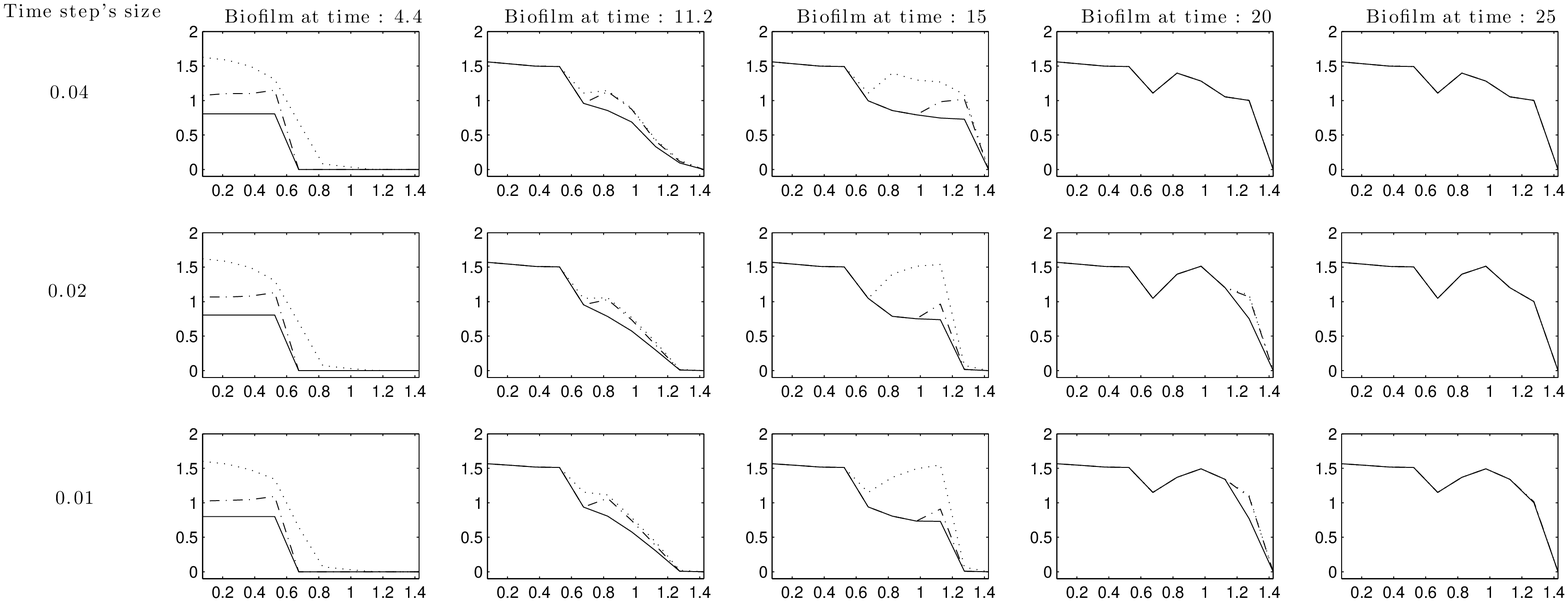}
\caption{\emph{Simulation 1 - 
Generation of two terraces with different time steps.}}
\label{figcompp}
\end{figure}

The goal of this first simulation is to show the ability of the model to generate
terraces and to study the stability of terrace generation with respect to time step.
We set the following values of the parameters 
\begin{gather} 
\begin{aligned}
\xic &= 0.1, & \tac &= 1,  & \bioz &=1,   &  \qcz &=0.2, &
  \gamt &=0.5,&  \gamd &= 0.9, \\
\etac &= 0.5, & A_w &= 1,     &  A_d &= 3.5 & 
   \kac &=2 ,& \alc &= 0.3, &\alc' &= 0.28,
\end{aligned}
\end{gather}
we take as coefficient $c(\hhc)$  involved  in swarmer cell velocity 
\begin{gather}
\label{defc} 
c(\hhc) = 0.2 \text{ if } \hhc<0.5 \text{ and } 0 \text{ if  } \hhc\geq 0.5,
\end{gather}
as position domain $[0,\posmax]= [0,1.5]$, and as initial functions 
\begin{align}
&\qcinit(\xc) = 0.1 \text{ if } \xc\leq 0.6 \text{ and } 0 \text{ otherwise,}
\label{sim1Q0}
\\
&\zec_0(\age,\xc)= 0 \text{ for all } \age \in [0,+\infty) \text{ and }\xc  \in [0,1.5] ,
\label{sim1Zeta0}
\\
&\rhc_0(\age,\bge,\xc)= 0 \text{ for all } \age \in [0,+\infty),\,
\bge \in [0,+\infty) \text{ and }\xc  \in [0,1.5].
\label{sim1Rho0}
\end{align}
The initial water concentration in the extra cellular matrix is set to $0.7$ 
and the initial water concentration in agar is set to 1.
The position step is set to  $\posst =0.15$ leading to a value of $I=11$  
and simulations are done with values of $\timest$ set to $0.04$, $0.02$
and $0.01$, with $\agest = \bgest = \timest$.
The results associated with those values of $\timest$ are given in figure \ref{figcompp} 
(line 1 gives the result for $\timest= 0.04$, line 2 gives the result for $\timest= 0.02$ and
line 3 gives the result for $\timest= 0.01$).
In each frame of this figure, the horizontal axis is the position axis and the vertical one gives the 
biomass densities.
The solid line is the vegetative cell density $\qc$, the 
dotted line represents the sum of the vegetative cell density and of 
the elongating cell biomass density $\qc + \biomm$ if the latter is 
non zero. The 
dashed line is the sum of the vegetative cell density,  of 
the elongating cell biomass density and of the swarmer cell biomass density
$\qc + \biomm + \biomn$, if the latter is non zero.

We see that for the three values of the time step, the dynamics of the colony and 
its final shape are similar. 
The more significant differences (see last column)  are located near $x=1.5$. 
Those differences are due to the influence of the boundary condition and
certainly also to a numerical drift linked  to the fact 
that this part of the colony is generated after a long term evolution.

Moreover, from a quantitative point of view,  if $\biom_i^n(\timest= 0.01)$ denotes the thickness of the 
colony at time $t_n$ and in position $x_i$ for the simulation made with $\timest= 0.01$, we have
\begin{gather}
\frac{\posst}{\posmax} \sqrt{\sum_{i=1}^{I}
\big( \biom_i^{25}(\timest= 0.02) -  \biom_i^{25}(\timest= 0.04)\big)^2 } = 0.0284,
\\
\frac{\posst}{\posmax} \sqrt{\sum_{i=1}^{I}
\big( \biom_i^{25}(\timest= 0.01) -  \biom_i^{25}(\timest= 0.02)\big)^2 } = 0.0181,
\end{gather}
suggesting  that when values of  time step has an order of magnitude of 0.01 or 0.02 we are
not far from the stability of the swarm colony dynamics with respect to time step.

\subsection{Simulation 2}
\begin{figure}
\centering 
\includegraphics[width=7.5cm]{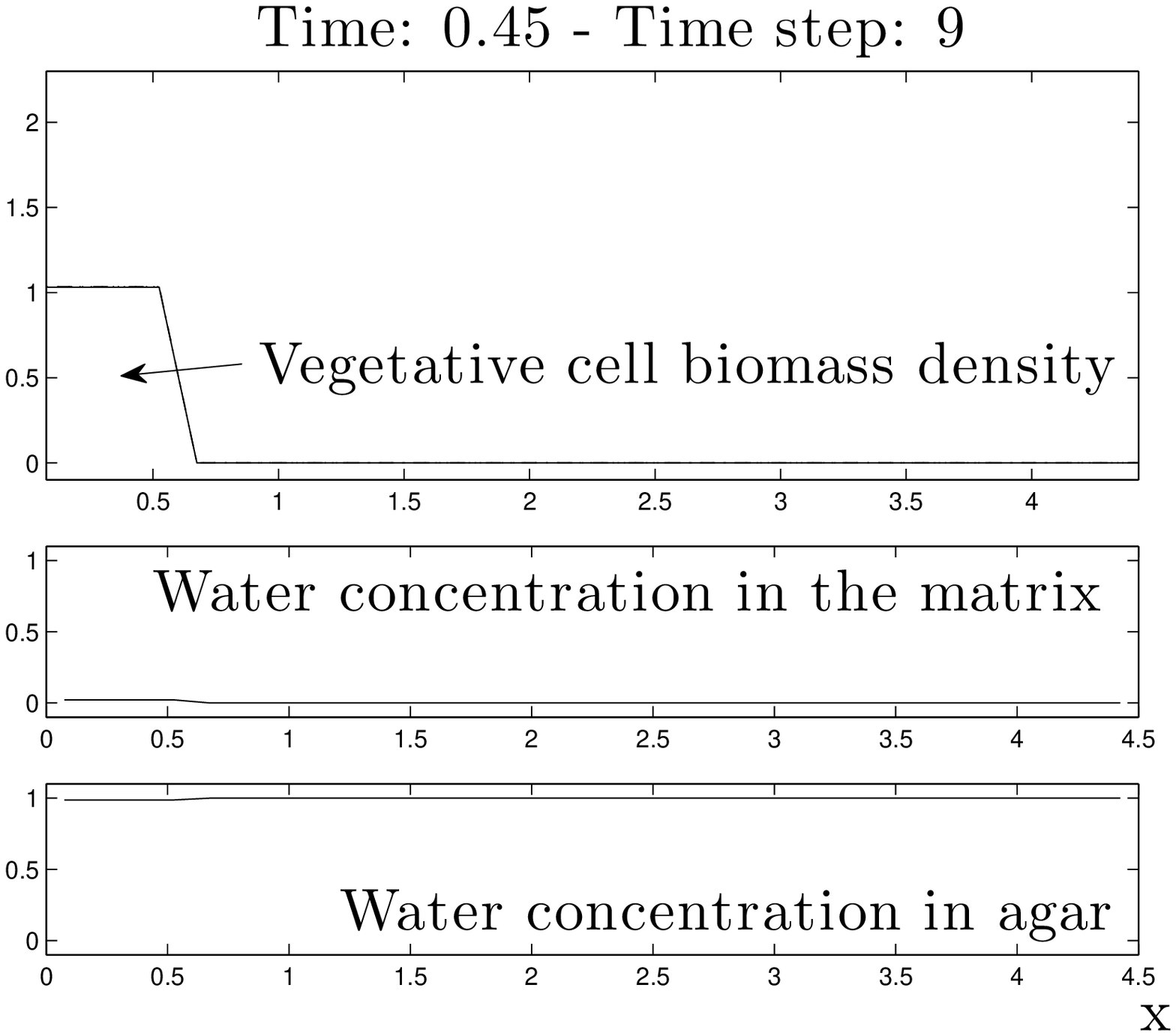}
\includegraphics[width=7.5cm]{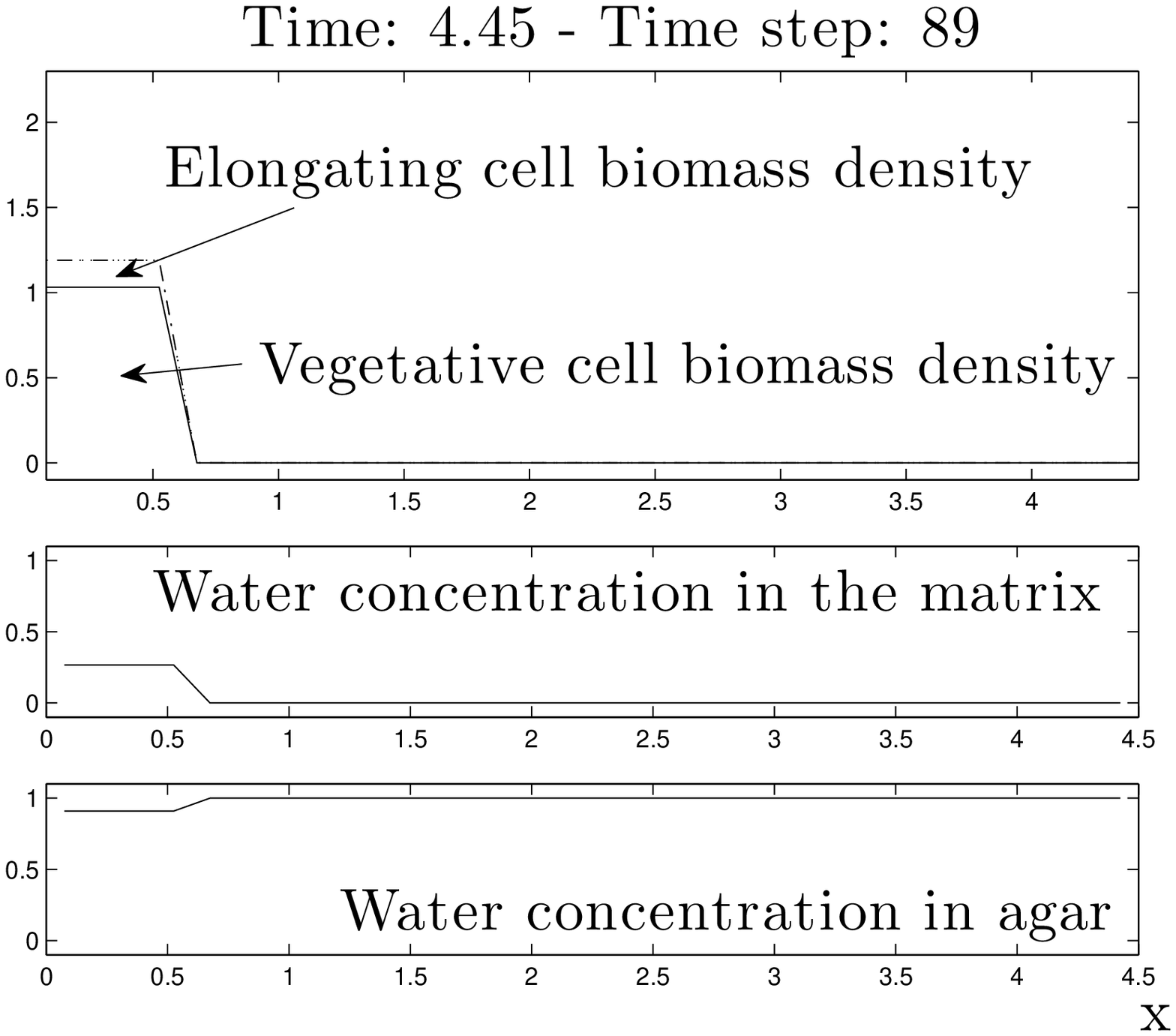}\\
\includegraphics[width=7.5cm]{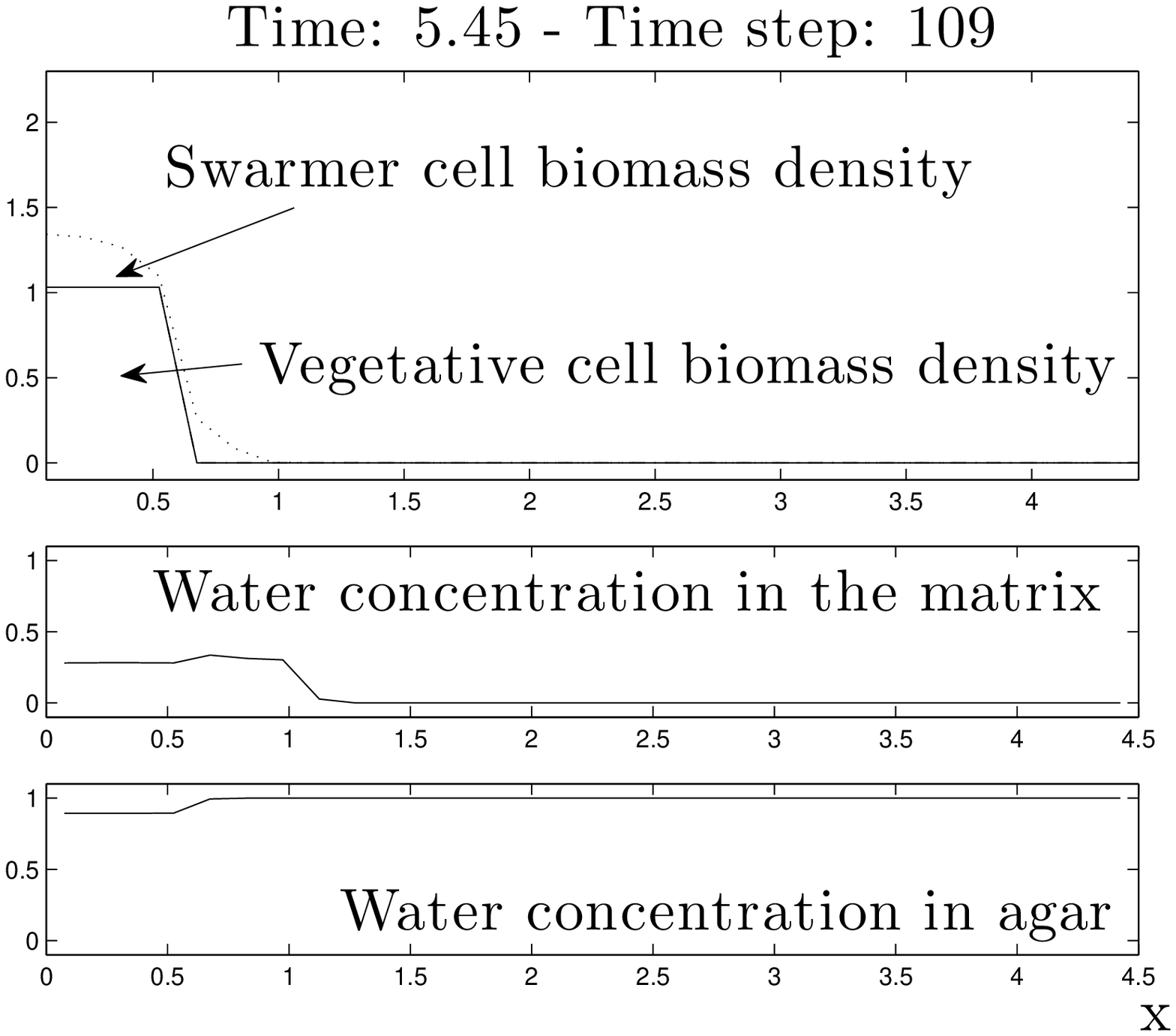}
\includegraphics[width=7.5cm]{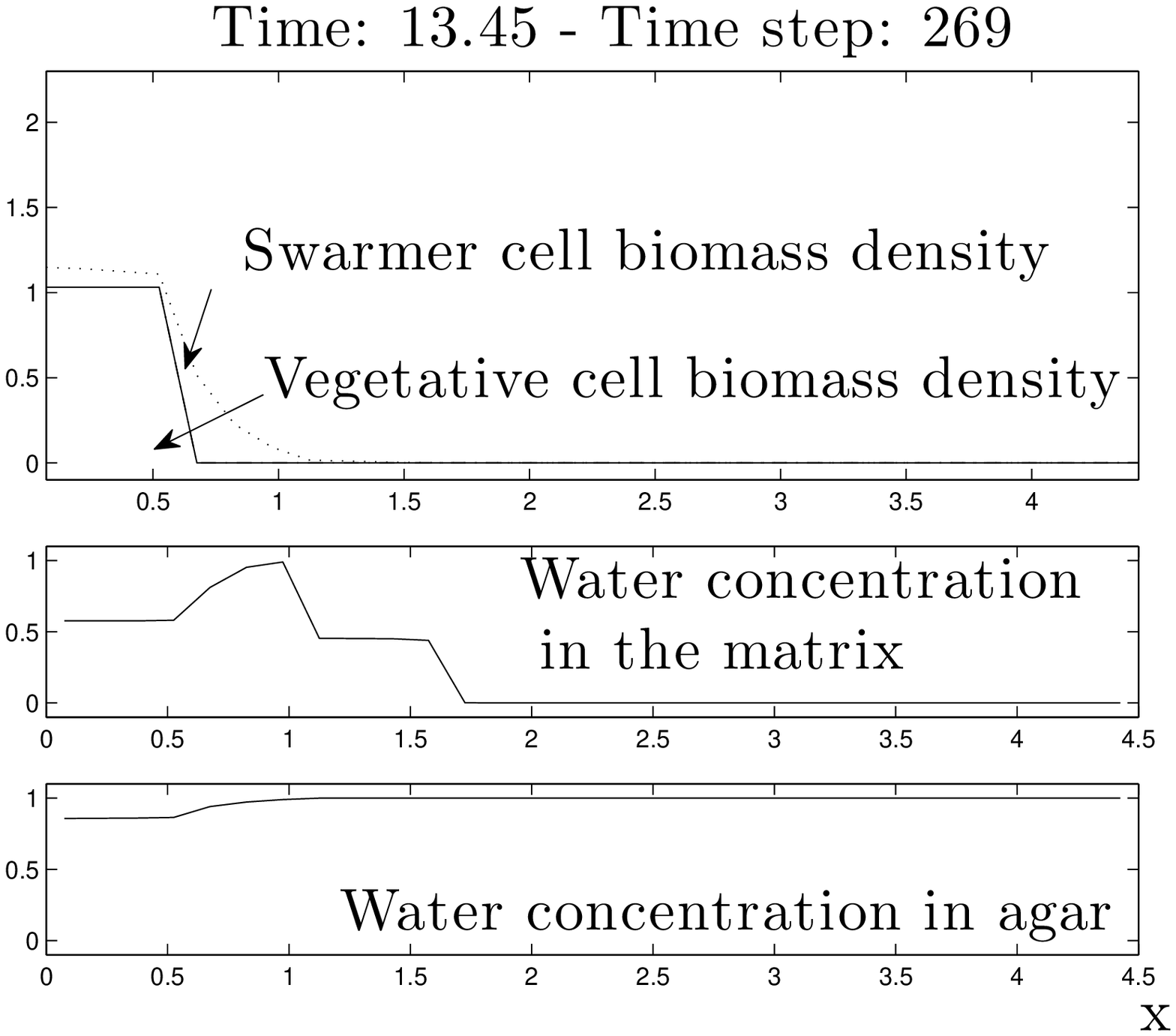}\\
\includegraphics[width=7.5cm]{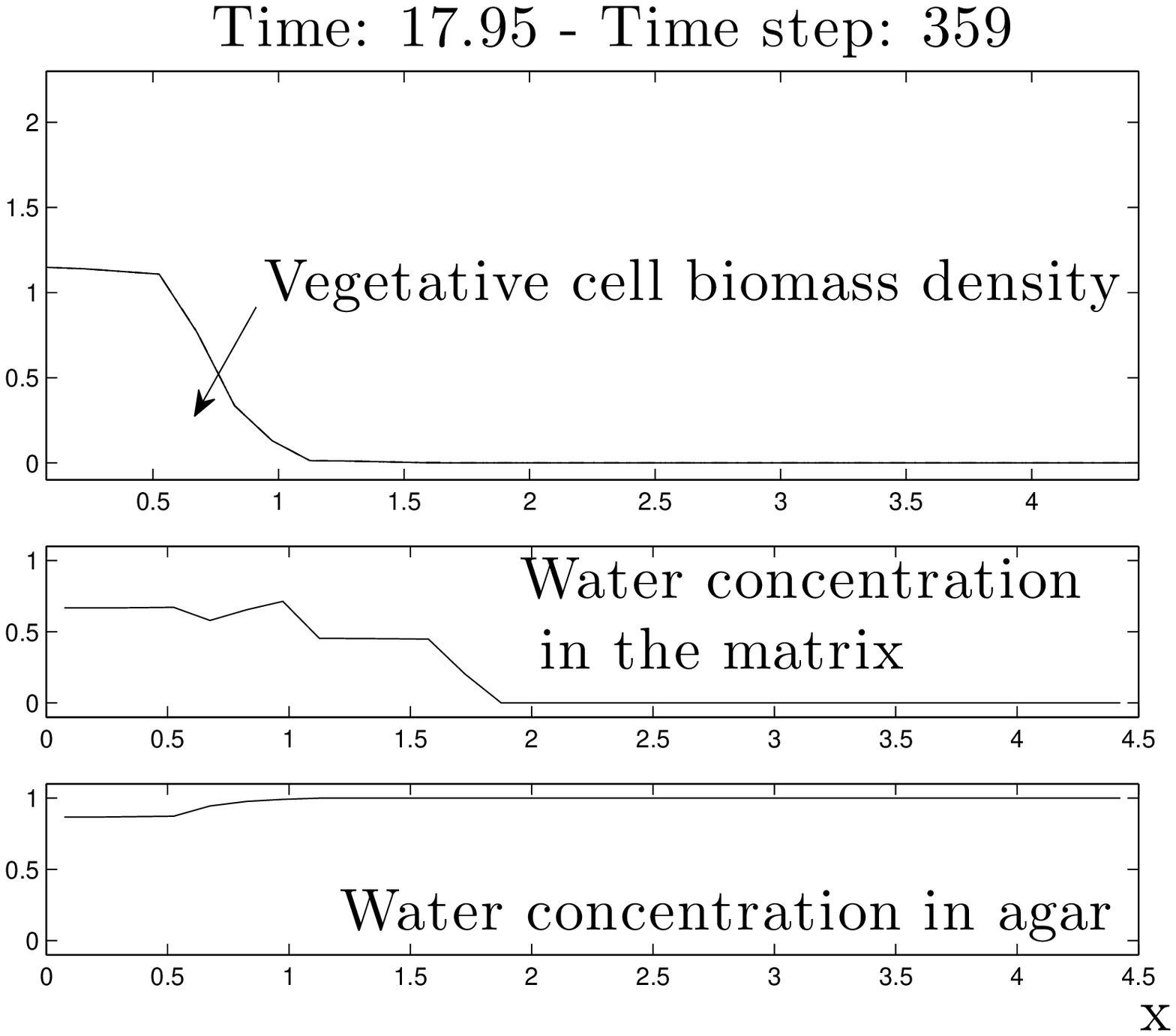}
\includegraphics[width=7.5cm]{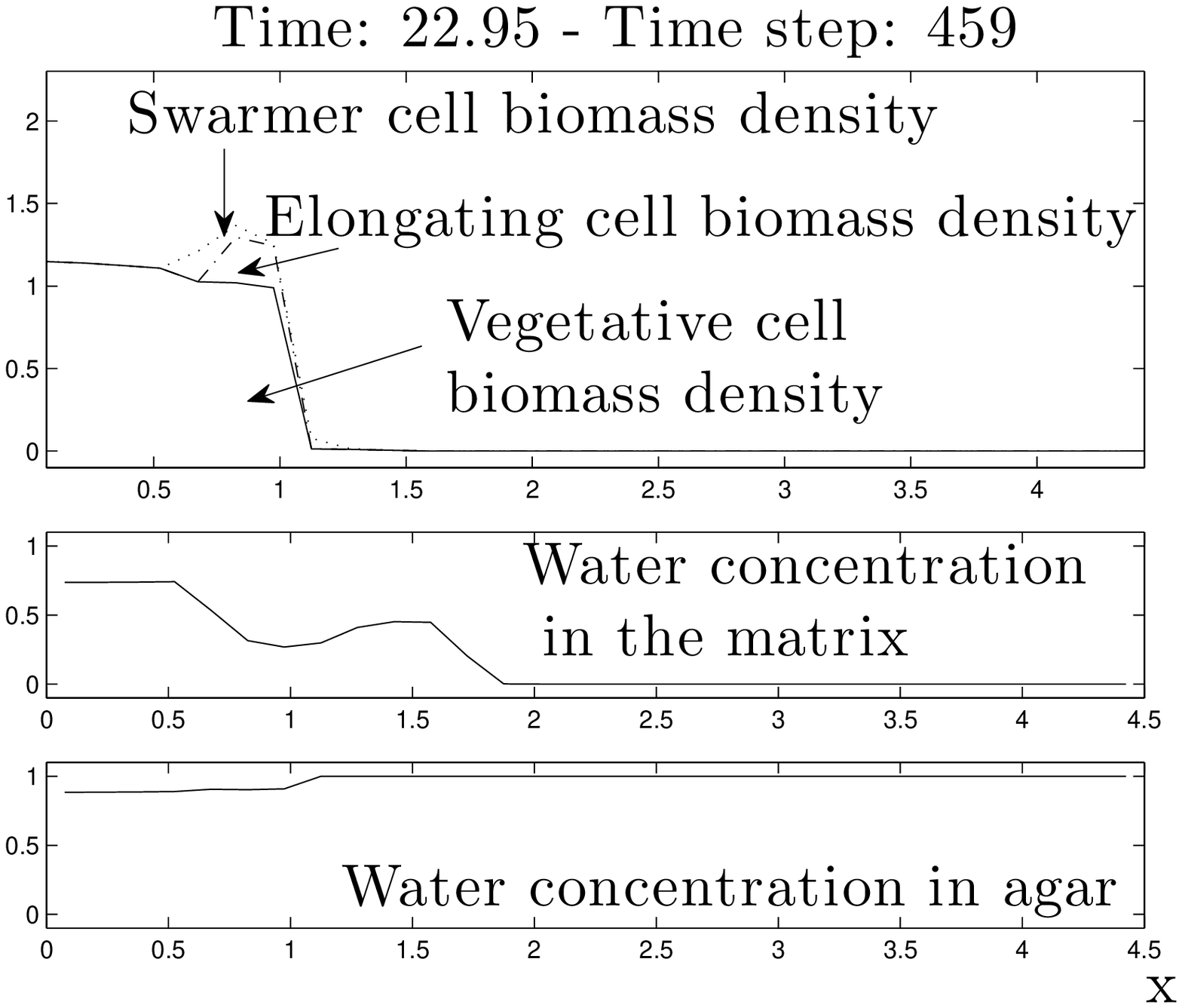}
\caption{\emph{Simulation 2 - Beginning of the swarm phenomenon:  
generation of the first terrace.}}
\label{figSim1A}
\end{figure}

The second simulation we present, consists in setting the following values of the
parameters
\begin{gather} 
\begin{aligned}
\xic &= 0.007, & \tac &= 1,  & \bioz &=1,   &  \qcz &=0.05, &
  \gamt &=0.03,&  \gamd &= 0.07, \\
\etac &= 0.3, & A_w &= 1,     &  A_d &= 6.3 & 
   \kac &=2.5,& \alc &= 0.02, &\alc' &= 0.0194,
\end{aligned}
\end{gather}
in taking as coefficient $c(\hhc)$  involved  in swarmer cell velocity the function
defined by (\ref{defc}), in taking as position domain $[0,\posmax]= [0,4.5]$,
and in taking as initial functions 
\begin{align}
&\qcinit(\xc) = 0.7 \text{ if } \xc\leq 0.6 \text{ and } 0 \text{ otherwise,}
\label{sim1Q0}
\\
&\zec_0(\age,\xc)= 0 \text{ for all } \age \in [0,+\infty) \text{ and }\xc  \in [0,4.5] ,
\label{sim1Zeta0}
\\
&\rhc_0(\age,\bge,\xc)= 0 \text{ for all } \age \in [0,+\infty),\,
\bge \in [0,+\infty) \text{ and }\xc  \in [0,4.5].
\label{sim1Rho0}
\end{align}
The initial water concentration in the extra cellular matrix is set to 0
and the initial water concentration in agar is set to 1.

Concerning the discretization parameters,  $\timest = \agest = \bgest =0.05$, 
$\posst =0.15$ leading to a value of $I=30$. The simulation is provided 
over 3000 time steps or, in other words, over the time interval $[0,150]$.

\:

\begin{figure}
\centering 
\includegraphics[width=7.5cm]{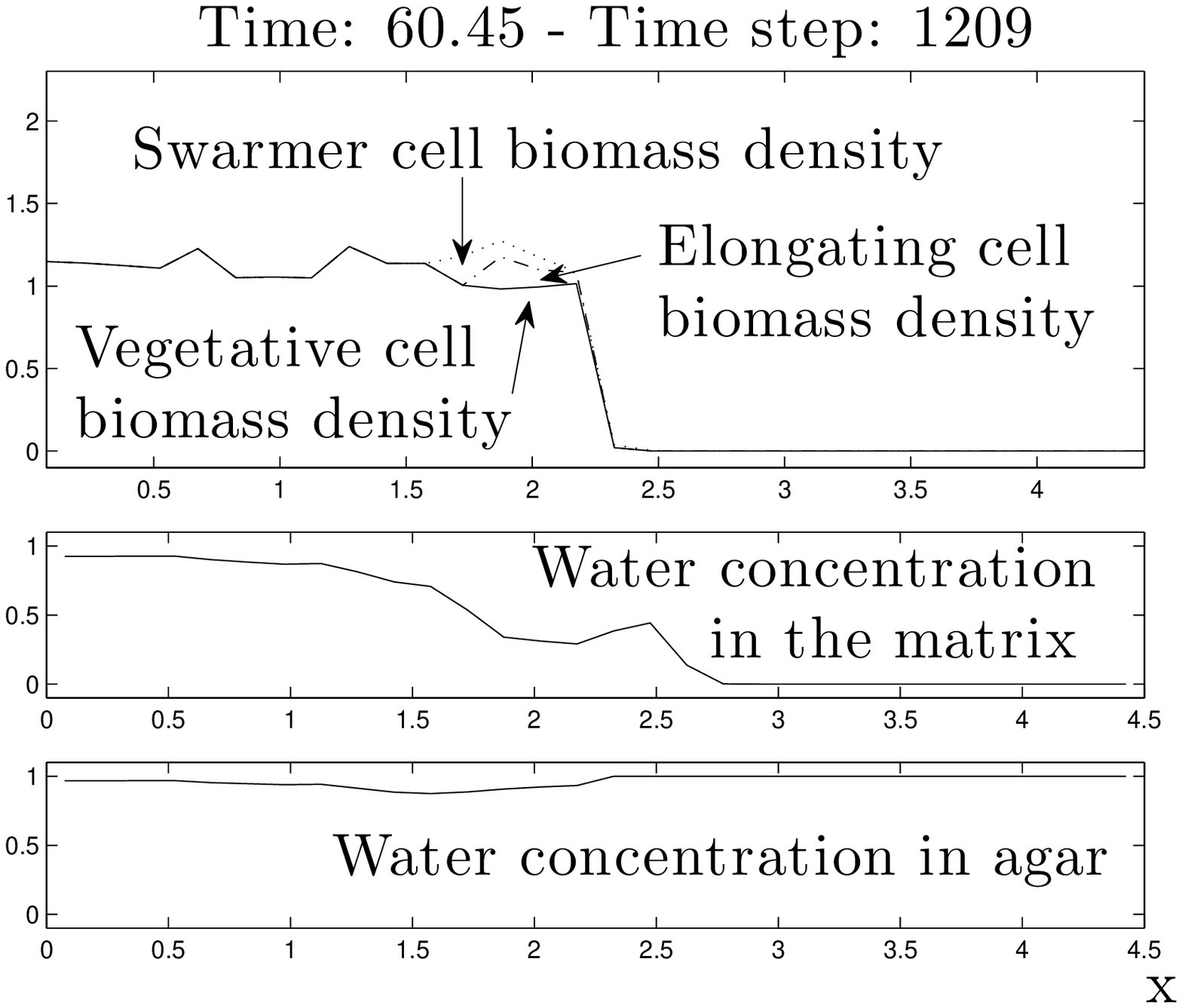}
\includegraphics[width=7.5cm]{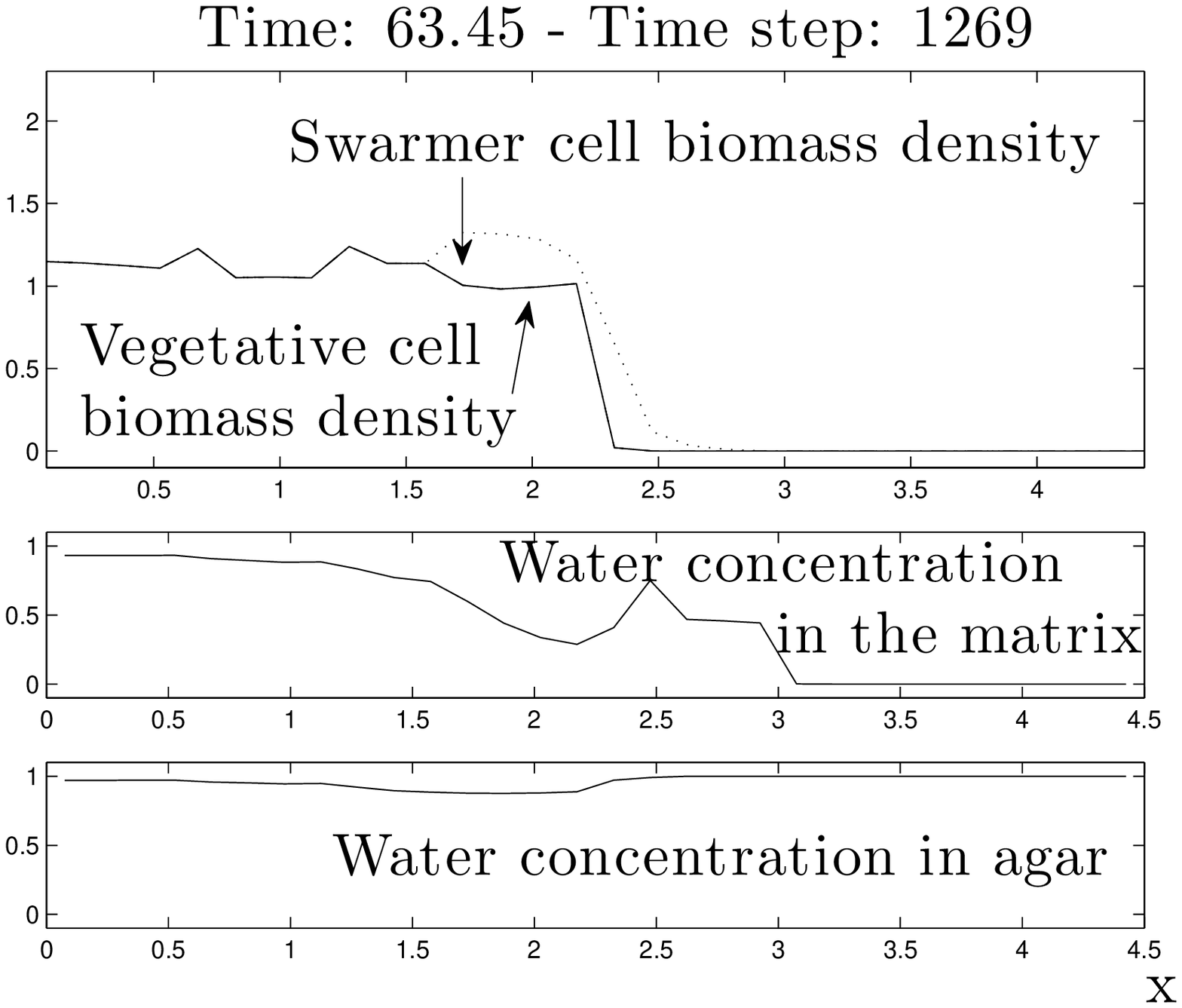}\\
\includegraphics[width=7.5cm]{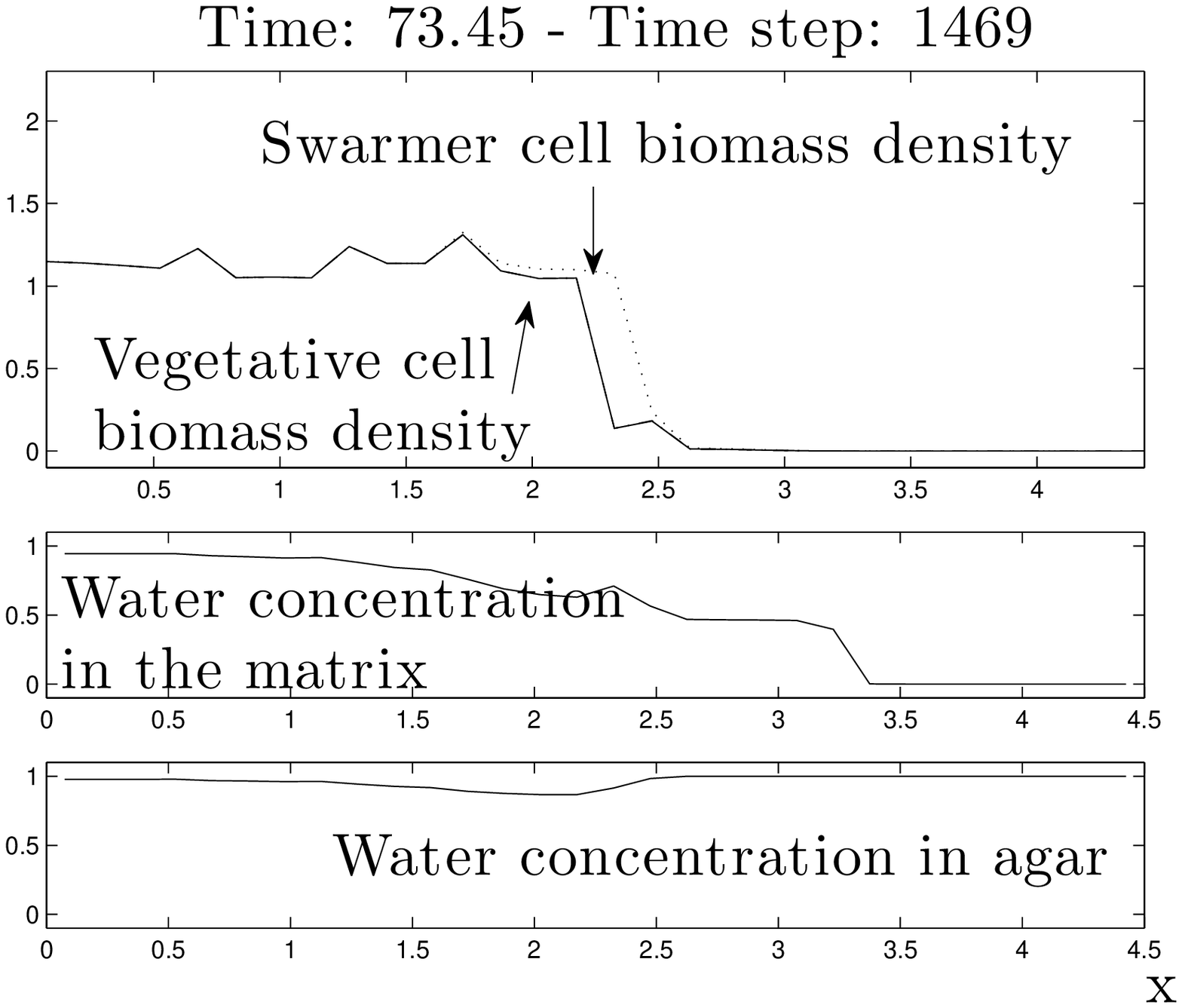}
\includegraphics[width=7.5cm]{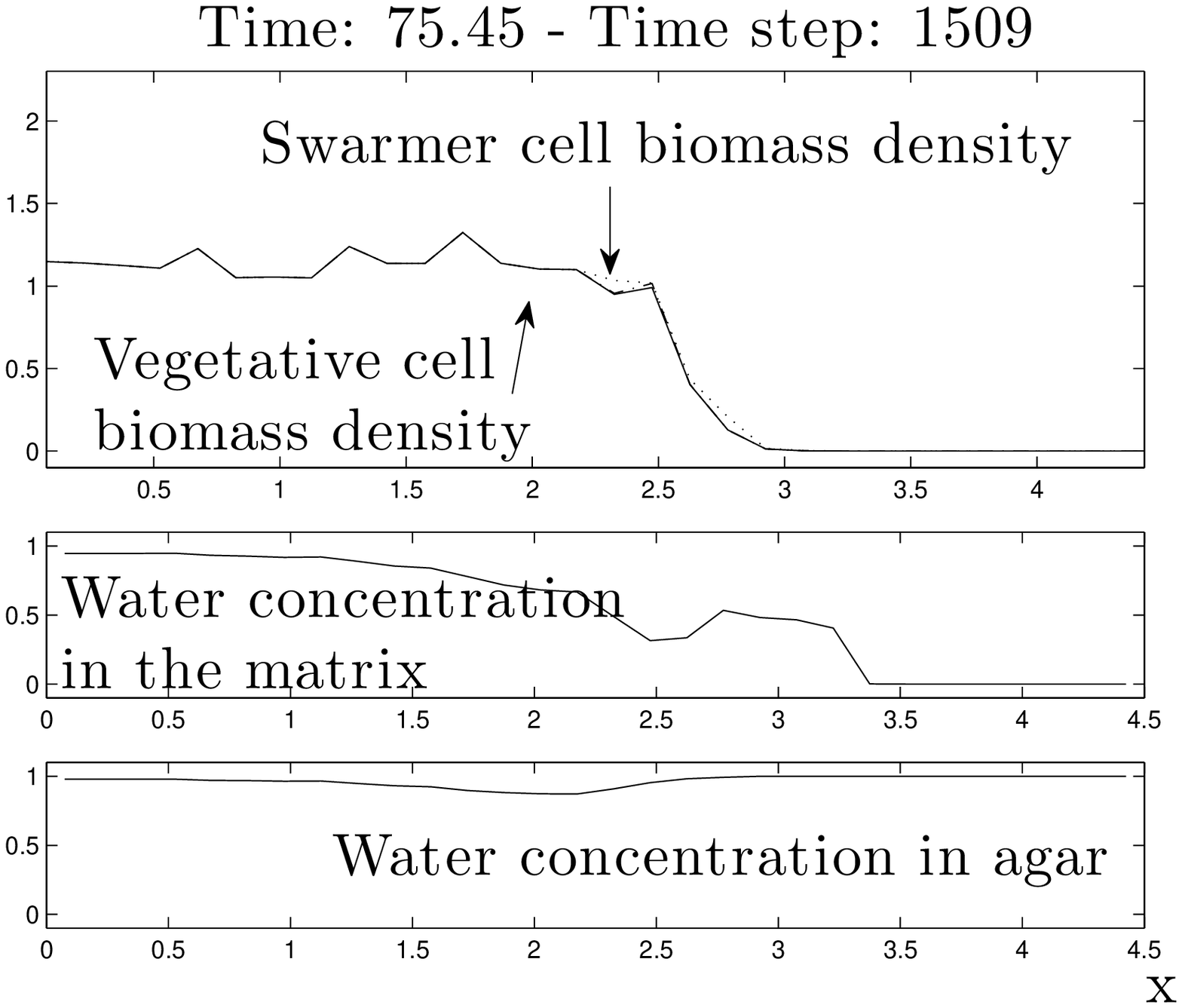}
\caption{\emph{Simulation 2 - 
Generation of a terrace in the middle of the swarm phenomenon.}}
\label{figSim1B}
\end{figure}

The beginning of the swarm phenomenon is summarized in Figure \ref{figSim1B}.
This figure is divided into six pictures which are in turn divided
into three frames.
The top frame shows the biomass densities at the time given just above it. 
The horizontal axis is the position axis and the vertical one gives the 
biomass densities.
The solid line is the vegetative cell density $\qc$, the 
dotted line represents the sum of the vegetative cell density and of 
the elongating cell biomass density $\qc + \biomm$ if the latter is 
non zero. The 
dashed line is the sum of the vegetative cell density,  of 
the elongating cell biomass density and of the swarmer cell biomass density
$\qc + \biomm + \biomn$, if the latter is non zero.
As a matter of fact, for instance, in the second picture the amount between 
the solid line and the dashed line has to be interpreted as the elongating cell
biomass density; and, in the third picture the amount between 
the solid line and the dotted line has to be interpreted as the swarmer cell
biomass density. The middle frame gives the water concentration $\hhc$, 
ranging between 0
and 1,  in the extra cellular matrix. The  bottom frame shows the water 
concentration $\ggc$ in agar.

The first picture shows the situation at time 0.45. 
From the initial time on, the vegetative cell density has increased
which has induced a  water consumption 
proportional to the vegetative cell density with the proportion factor
being $\alc$. As a consequence, at that time, the extra cellular matrix is very dry.
Hence, the elongating cells which are created have a long existence
duration expectancy. At the given time, the thickness just reaches 1. 
Hence the cell division
stops and the only way for the colony thickness to grow  is to increase 
the elongating cell biomass density. This is made possible by the elongation
of each elongating cell, the mass of which then grows,
and by the fact that elongating cells have a large life expectancy.

The second picture shows the situation at time 4.45. Since time 0.45,
the elongating cell biomass density, and the colony thickness
have been increasing. The water concentration in agar has decreased little
due to exchange with the extra cellular matrix.
In the extra cellular matrix, the water consumed by vegetative cells
is proportional to the vegetative cell density with the proportion factor 
being $\alc -\alc'$. The water consumed by elongating cells
is proportional to the elongating cell biomass density with the proportion factor 
being $\alc$. As a result of the water consumption and  exchanges, 
at this particular time, the extra cellular matrix is
little hydrated.
Soon after time 4.45, elongating cells become swarmer cells and, as the 
extra cellular matrix water concentration is low, the first swarm step begins.

The situation at time 5.45 (see third picture) is the following: 
the first swarm step is taking place and then swarmer cells move towards
the positive $\xc-$axis. At the same time, the extra cellular
matrix water concentration near $\xc=0.8$, \emph{i.e.} near the colony edge,
is increasing faster than near the colony center. 
As a consequence, soon after time 5.45 the extra cellular 
matrix water concentration near the colony edge goes beyond the threshold
$H=0.5$ stopping swarmer cells, and then jamming the first swarm step.

Between time 5.45 and time 13.45 (see fourth picture), the swarm colony shape has 
not evolved a great deal.
This period of time has to be interpreted as
the first consolidation phase. 
The space distribution of the cell size, in accordance with  
Matsuyama \emph{et al.} \cite{MatTaItWaMat}, shows a large proportion of 
long cells in the colony front and a smaller proportion in the colony interior.
On the other hand, the extra cellular matrix 
water concentration has grown. In particular the  threshold $H=0.5$ 
has been overtaken in the colony center, which will prevent the swarm from 
going towards the colony center while the next swarm step begins.

At time 17.45 the first consolidation phase ends
(see fifth picture). Since time 13.45,
every swarmer cell has been de-differentiating to produce new vegetative cells,
making the proportion of short cells to increase gradually, always in accordance with  
Matsuyama \emph{et al.} \cite{MatTaItWaMat}.
 As the total biomass density (or the thickness), near the colony
edge, is low the cell division process may happen giving rise to an
increase of vegetative cell density.
At the same time, as the water consumption is important (proportional
to the vegetative cell density with proportion factor $\alc$), the 
water concentration of the extra cellular matrix decreases fast.

At time 22.95 (see sixth picture), near $\xc=1$,  the extra cellular matrix 
water concentration is low. Hence the elongating cell 
existence duration expectancy is long. This leads to an increase in thickness,
despite that the thickness is larger than 1.
Elongating cells turn into swarmer cells which cannot go
towards the colony center because of the high extra cellular matrix
water concentration in this region. Hence they go towards
positive $\xc-$axis, beginning the second swarm step.

\:

Figure \ref{figSim1B} depicts the generation of a terrace in the middle of the
swarm phenomenon.

The situation in the first picture (time 60.65) is very similar to the one in the
sixth picture of Figure \ref{figSim1A}: thickness is larger than 1,
thickness increases due to elongating cells, 
elongating cells turn into swarmer cells,
relative dryness of the 
extra cellular matrix near the colony edge and relative wetness of the 
extra cellular matrix while going towards the colony center which prevent 
swarming in that direction.

The second picture (time 63.45) shows the situation at the end of 
the swarm step, \emph{i.e.} at the beginning of the consolidation phase. 
The characteristics of this moment are: 
wetness of the extra cellular matrix near the colony edge and in 
the colony center, relative dryness with increasing water concentration
of the extra cellular matrix between the colony edge and the
the colony center.

The colony shape has not evolved much between time 63.45 and 
time 73.45 (third picture) but the extra cellular matrix water
concentration has increased.
At time 73.45 the de-differentiation of swarmer cells is just
beginning and producing vegetative cells, the population of which is
growing by cell division.  
Here again, in accordance with  Matsuyama \emph{et al.} \cite{MatTaItWaMat},
a large proportion of long cells is present in the colony front and gradually this 
proportion decreases with time.  

The water consumption is linked to the increase yield of a relatively dry
extra cellular matrix near the colony edge, preparing the next
swarm step beginning at time 75.45 (see forth picture).

In the third and forth pictures we can see that a remaining hump around
$x=1.7$ has been created during the time period covered by the pictures
of Figure \ref{figSim1B}. This hump constitutes a terrace.
\begin{figure}
\centering 
\includegraphics[width=7.5cm]{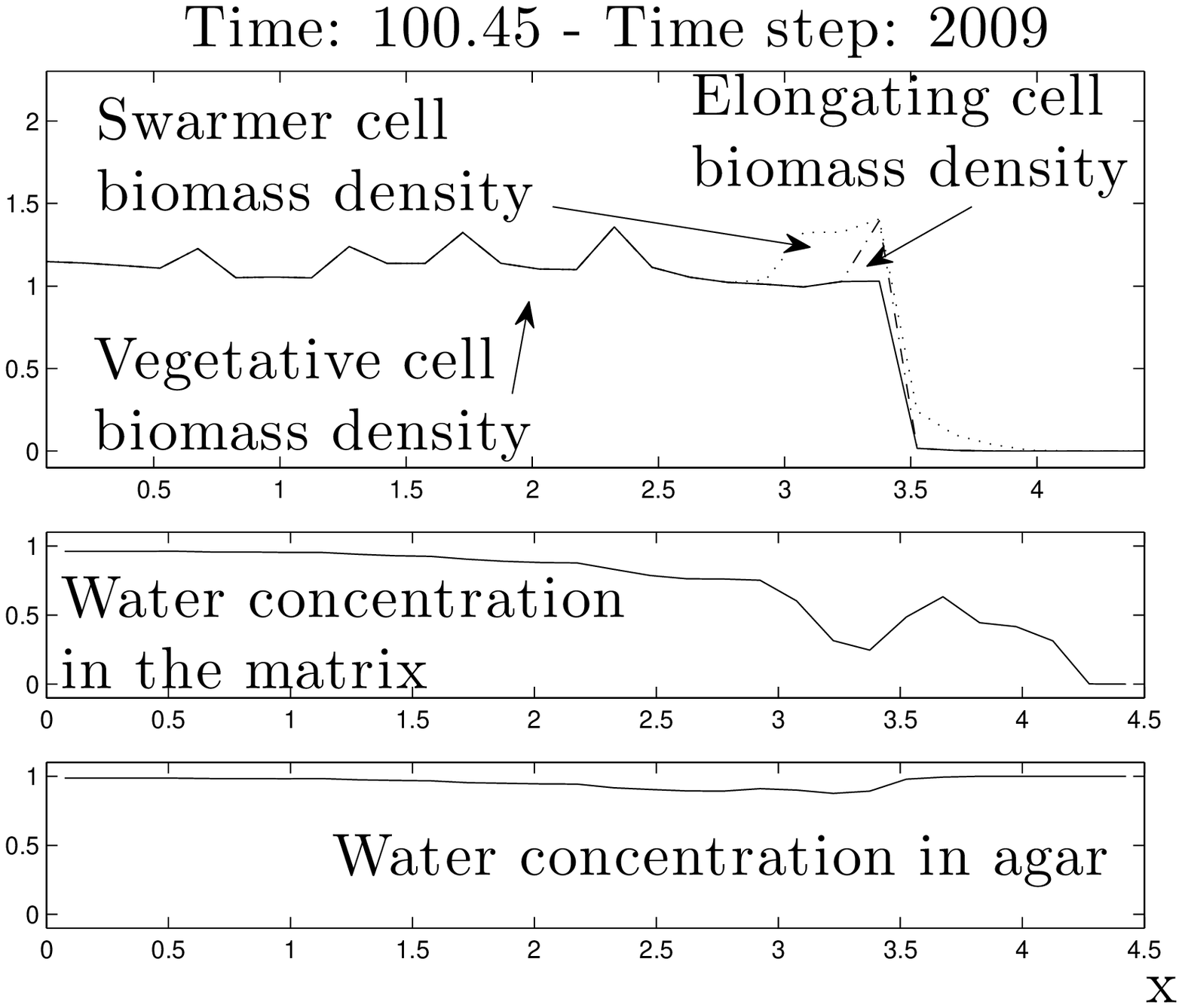}
\includegraphics[width=7.5cm]{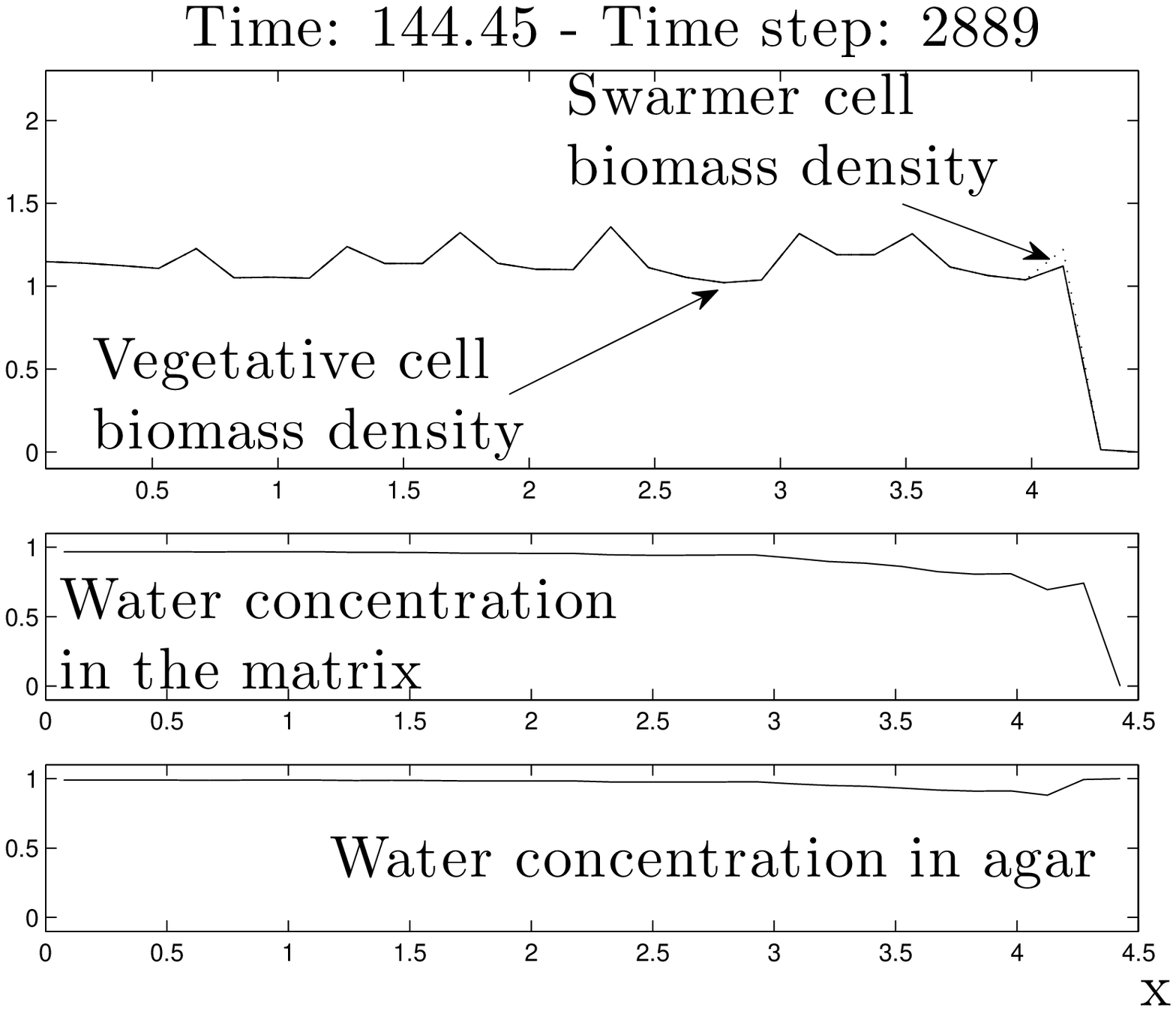}
\caption{\emph{Simulation 2 - End of the colonization. }}
\label{figSim1C}
\end{figure}

\:

This colonization goes on until the whole position domain is conquered.
Figure \ref{figSim1C} shows the situation at time 100.45 and 
close to the end of the colonization process.
The relative regularity of the terraces can be observed in these two pictures.
We can also notice that the water concentration in agar and in the matrix
slowly goes to 1. 
\begin{figure}
\centering 
\includegraphics[width=15cm,height=6cm]{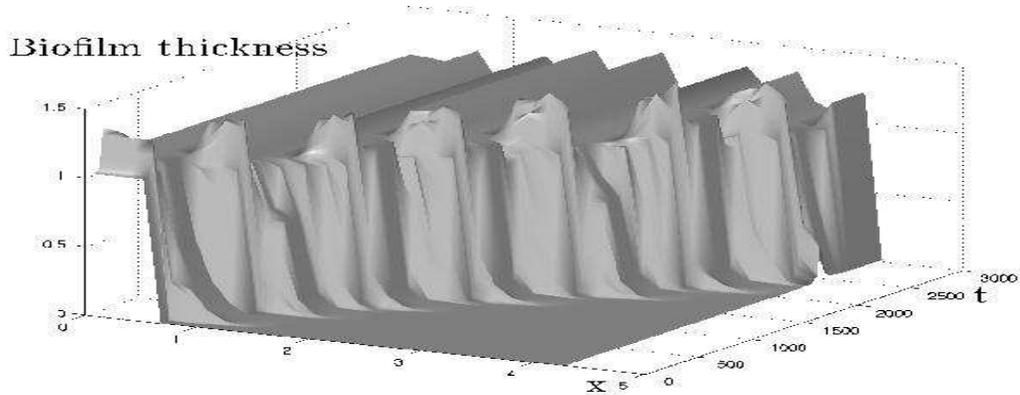}
\caption{\emph{Simulation 2 - Summarize of the swarm phenomenon.}}
\label{figSim1CC}
\end{figure}

\:

Figure \ref{figSim1CC} summarizes the swarm phenomenon over the
time period $[0,150]$ \emph{i.e.} from 0 to the $3000^\text{th}$
time step.
The horizontal axis going backward shows the time (in time step units)
and the horizontal axis going to the right shows the position.
The vertical axis represents the total biomass or thickness $\qc + \biomm + \biomn$.
This figure illustrates the relative regularity of the phenomenon,
in time and position.

Moreover, it also allows one to visualize how the terraces are formed: 
a large hump is generated. Then it drains towards the colony edge
but only partially, giving rise to a remaining hump which makes a 
terrace.

\:

\begin{figure}
\centering 
\includegraphics[width=11cm,height=7cm]{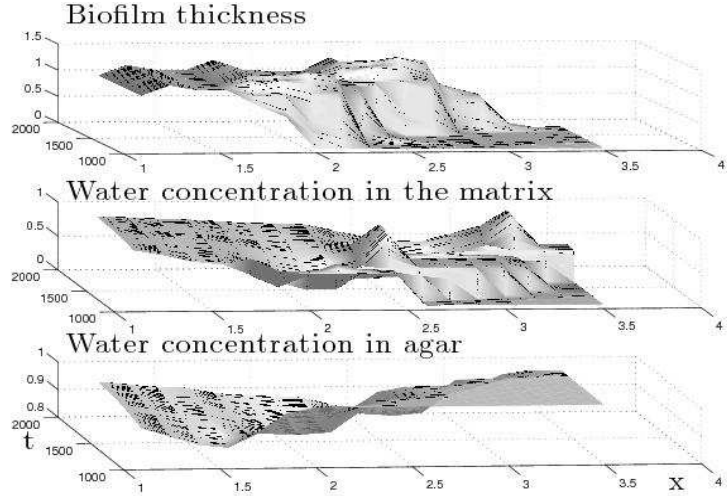}
\caption{\emph{Simulation 2 - Visualization of the relative positions of terraces, 
matrix water concentration maxima and agar water concentration minima.}}
\label{figSim1D}
\end{figure}

Finally, the first frame of Figure\ref{figSim1D} shows 
the same object as in Figure \ref{figSim1CC} but restrained
to time interval $[1000\timest, 2000\timest] = [50, 100]$.
The second frame shows the water concentration in the extra cellular
matrix and the third one shows the water concentration in agar.

This figure shows that,  during the terrace formation,
in front of the hump, the water concentration in the extra cellular
matrix is high; and that, at the back of the hump, the water 
concentration in agar is low.

\subsection{Simulation 3}
The third simulation we present is carried out with the same
position domain as before,
with:
\begin{gather} 
\begin{aligned}
\xic &= 0.008, & \tac &= 1,  & \bioz &=1,   &  \qcz &=0.05,& 
  \gamt &= 0.37, & \gamd &= 0.13 , \\
\etac &= 0.3, & A_w &= 1.2,     &  A_d &= 5.5, & \kac &=2.5, & 
   \alc &= 0.42, & \alc' &= 0.41,
\end{aligned}\\
\qcinit(\xc) = 0.2 \text{ if } \xc\leq 0.6 \text{ and } 0 \text{ otherwise,}
\label{sim1Q0}
\end{gather}
function $\zec_0$ given by (\ref{sim1Zeta0}), $\rhc_0$ by (\ref{sim1Rho0}) 
and $c$  by (\ref{defc}).
The initial water concentration in the extra cellular matrix is set to 0.8
and the initial water concentration in agar is set to 1.

The steps are  $\timest = \agest = \bgest =0.05$, 
$\posst =0.15$ leading to a value of $I=30$. The simulation is provided 
over 4400 time steps or over the time interval $[0,220]$.

Some parameters are quite different from the previous simulation. But the
main significant difference is that here   $\gamt$ is larger than $\gamd$.
\begin{figure}
\centering 
\includegraphics[width=7.5cm]{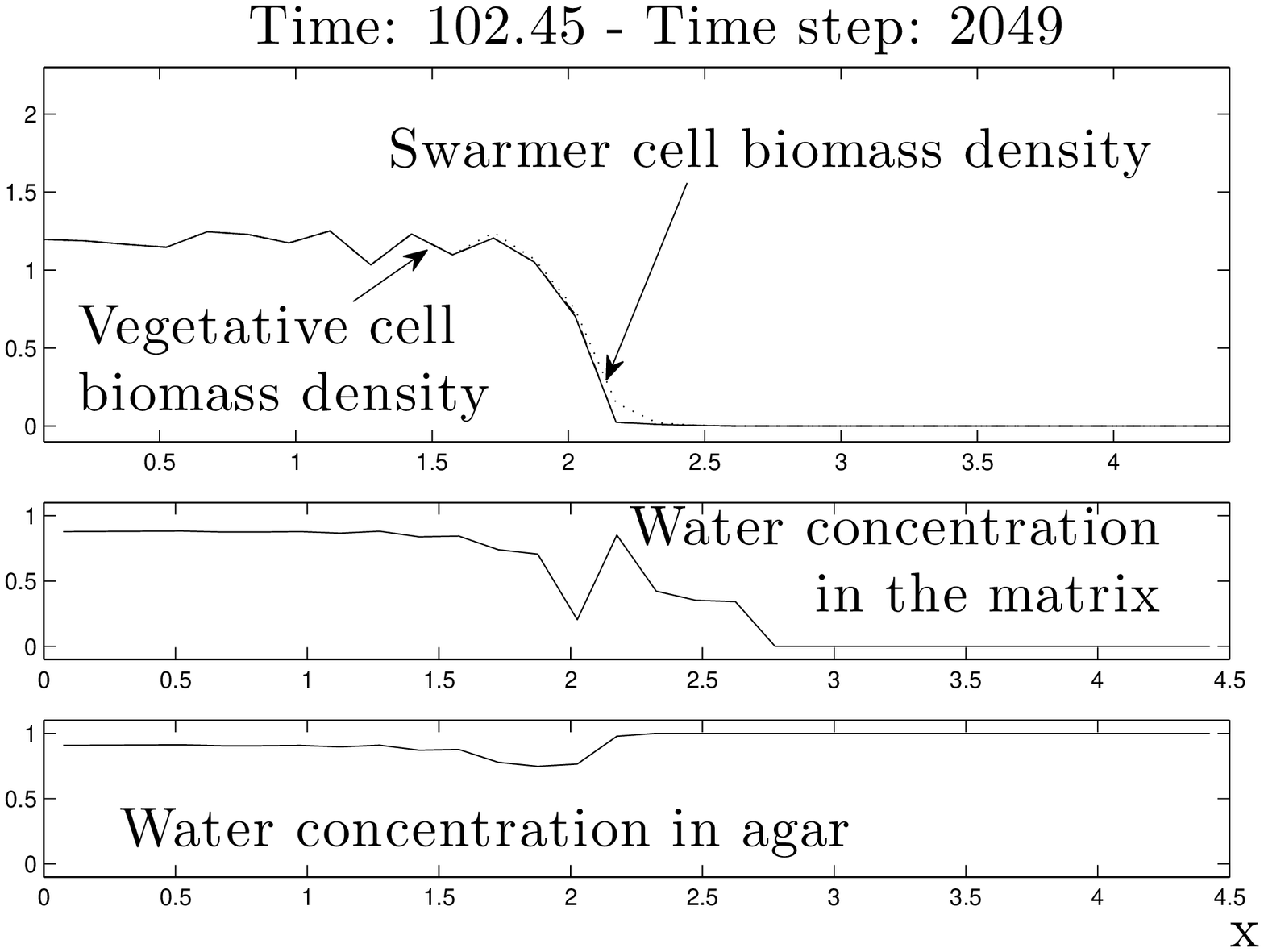}
\includegraphics[width=7.5cm]{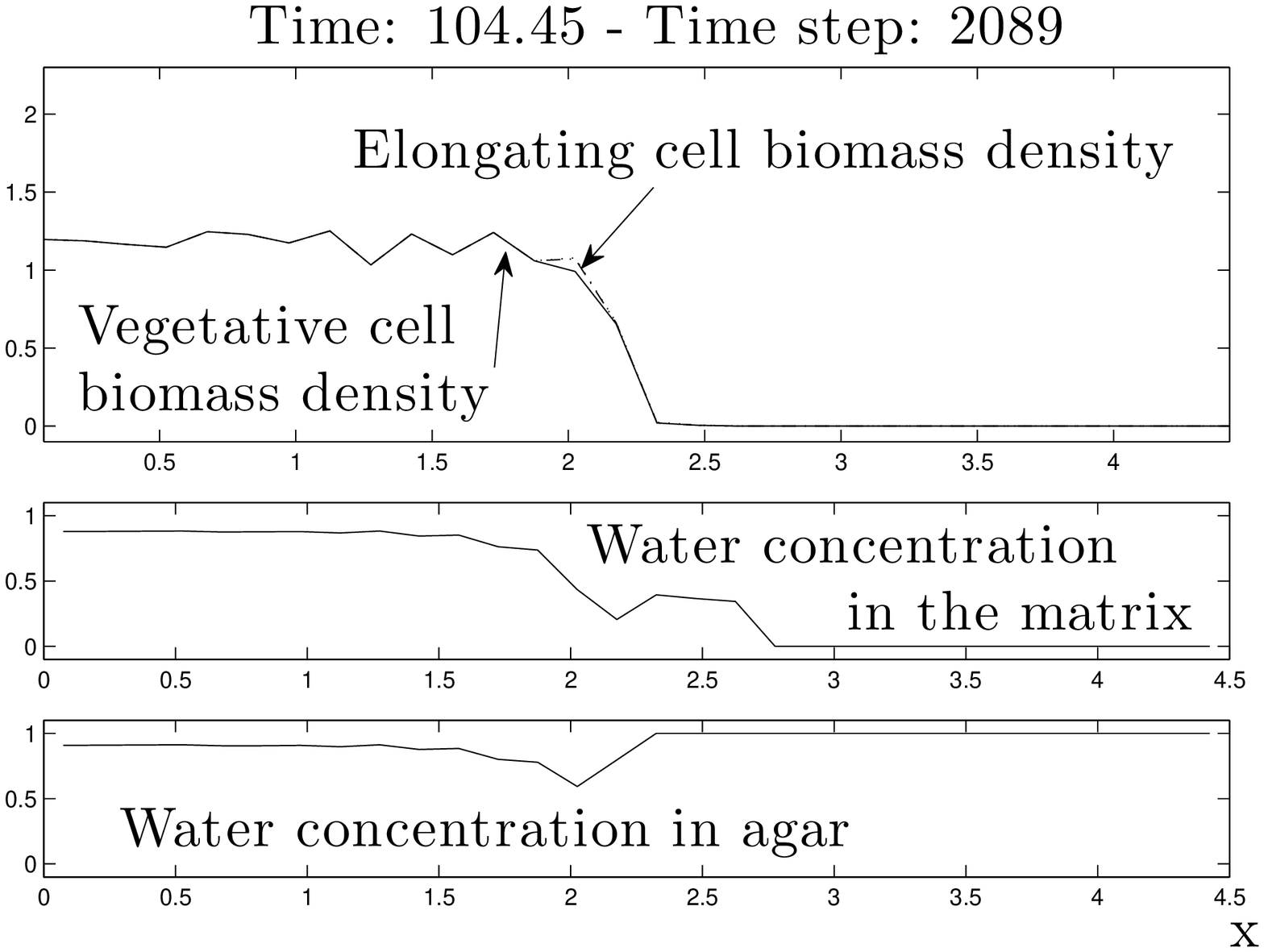}\\
\includegraphics[width=7.5cm]{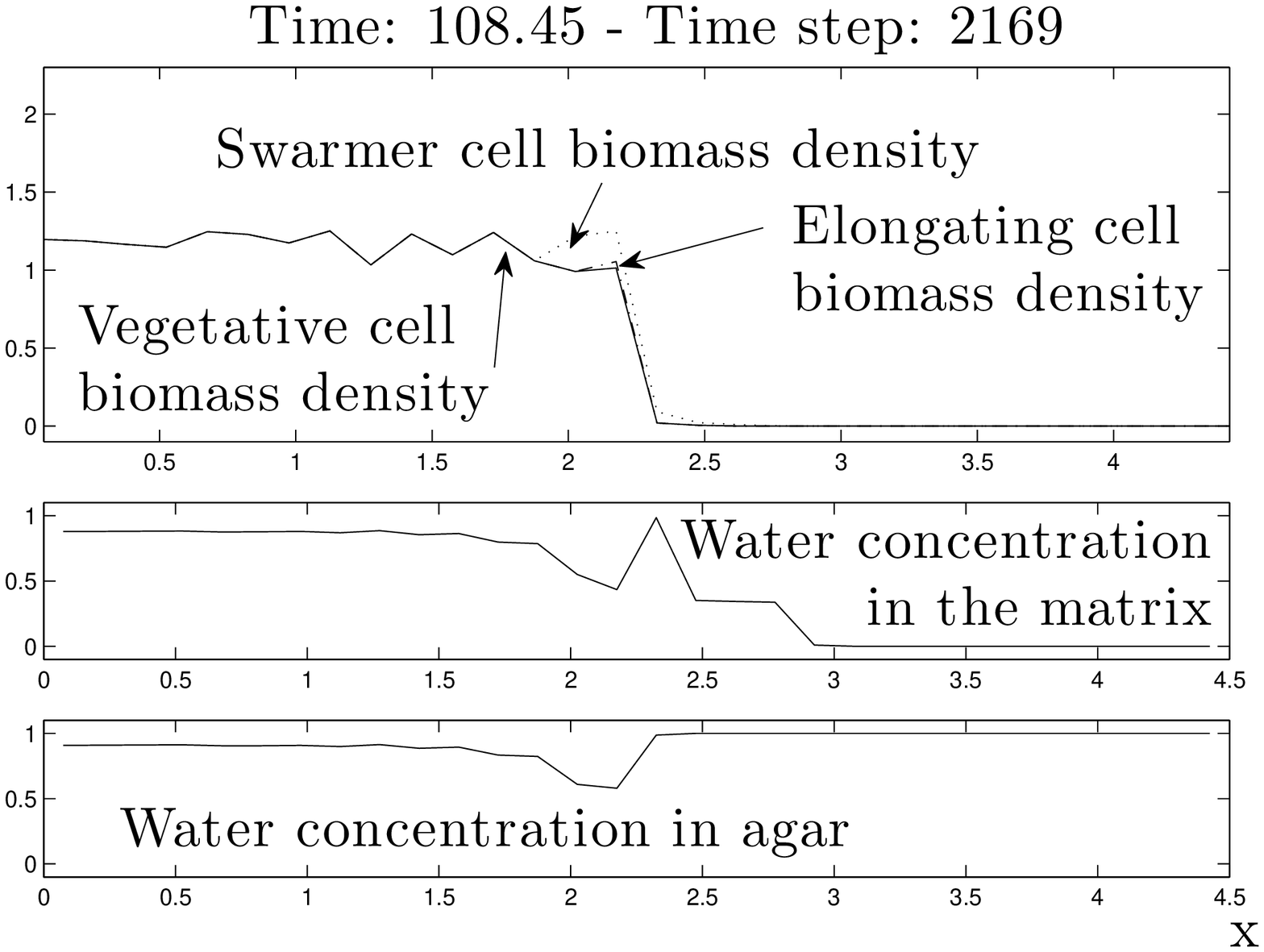}
\includegraphics[width=7.5cm]{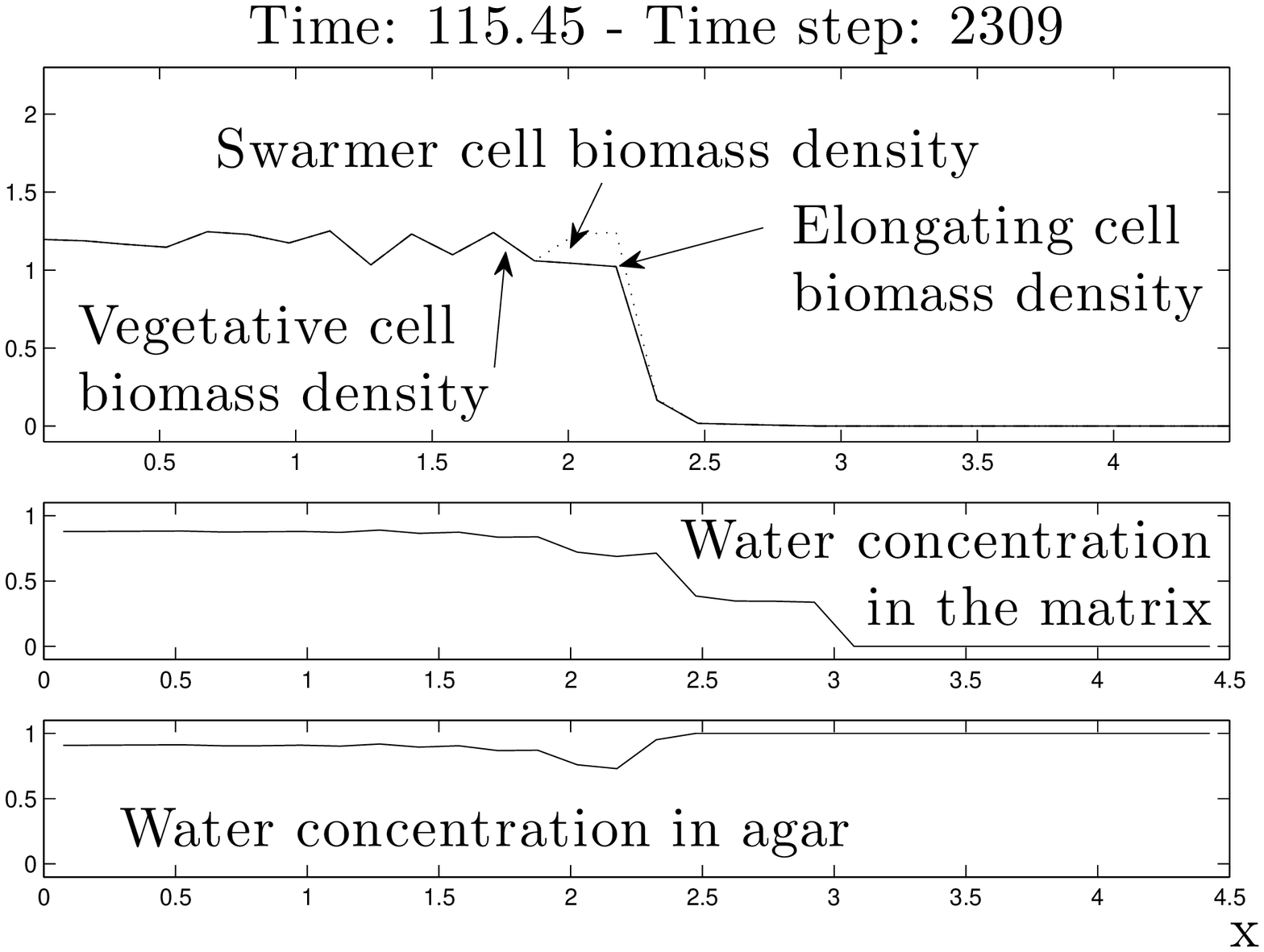}\\
\includegraphics[width=7.5cm]{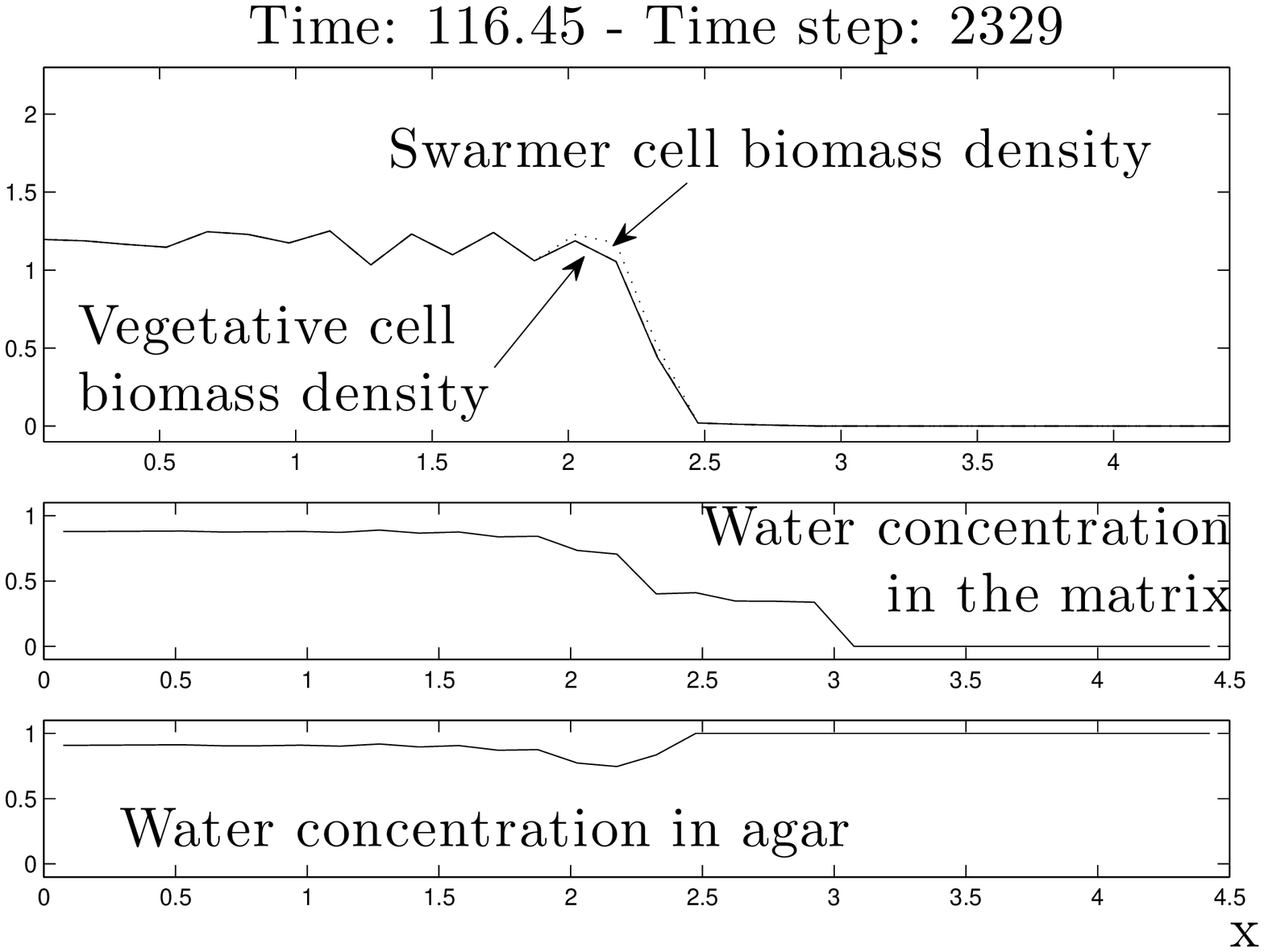}
\includegraphics[width=7.5cm]{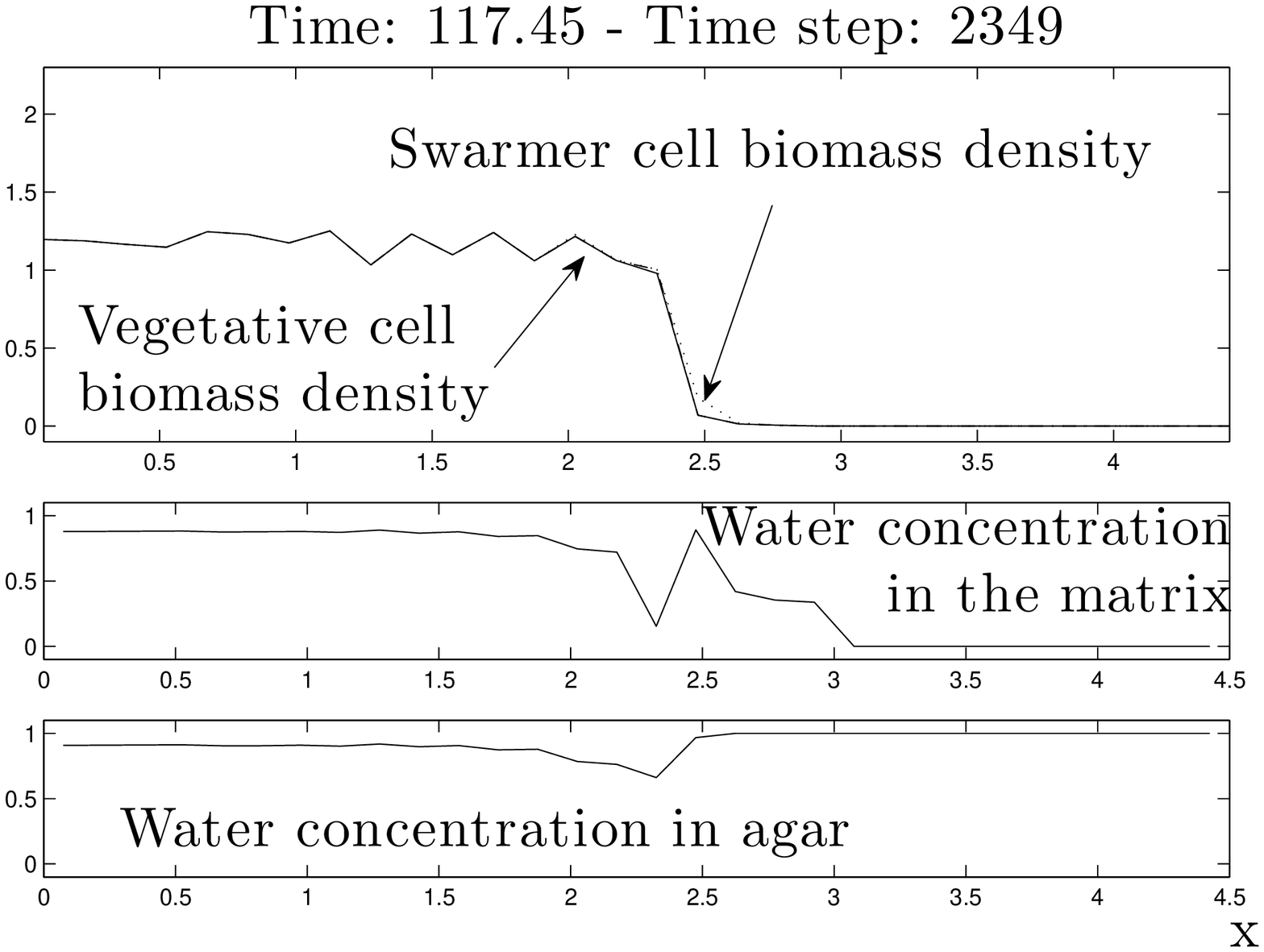}
\caption{\emph{Simulation 3 - Generation of a terrace.}}
\label{figSim2A}
\end{figure}
This has the following consequence, as is visible in Figure \ref{figSim2A} 
which shows a terrace generation:
the water which is contained in agar is used more.
Apart from this difference, the way a terrace is created is very similar 
to simulation 1.

The first picture shows the situation at a moment when swarmer cells are
de-differentiating and producing vegetative cells which are growing
by cellular division. This is drying the extra cellular matrix. 

The second picture shows a colony thickness, near the colony edge,
larger than 1. Then
thickness increases
 thanks to elongating cells which are 
turning into swarmer cells.
The extra cellular matrix near the colony edge is dry, while towards the colony center
it is wet.

The third picture  shows the situation at the beginning of the consolidation phase. 
At this moment, the extra cellular matrix near the colony edge and in 
the colony center is wet. Between the colony edge and the
the colony center the water concentration
of the extra cellular matrix increases. 

The swarm colony shape does not evolve much until the moment of the 
third picture but the extra cellular matrix water concentration increases.

The fifth picture shows the moment when the de-differentiation begins which
will soon lead to the situation in the sixth picture.

The situation in the sixth picture is similar to the one in the first picture
but with a colony edge which has moved towards the right. 
\begin{figure}
\centering 
\includegraphics[width=7.5cm]{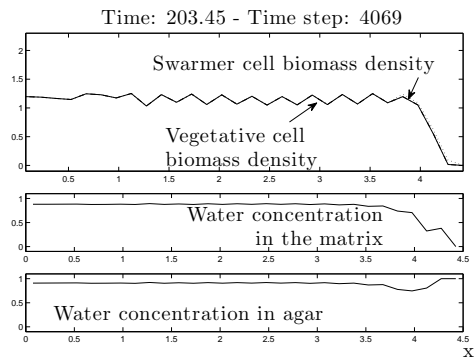}
\caption{\emph{Simulation 3 - End of the colonization.}}
\label{figSim2B}
\end{figure}

\begin{figure}
\centering 
\includegraphics[width=15cm,height=6cm]{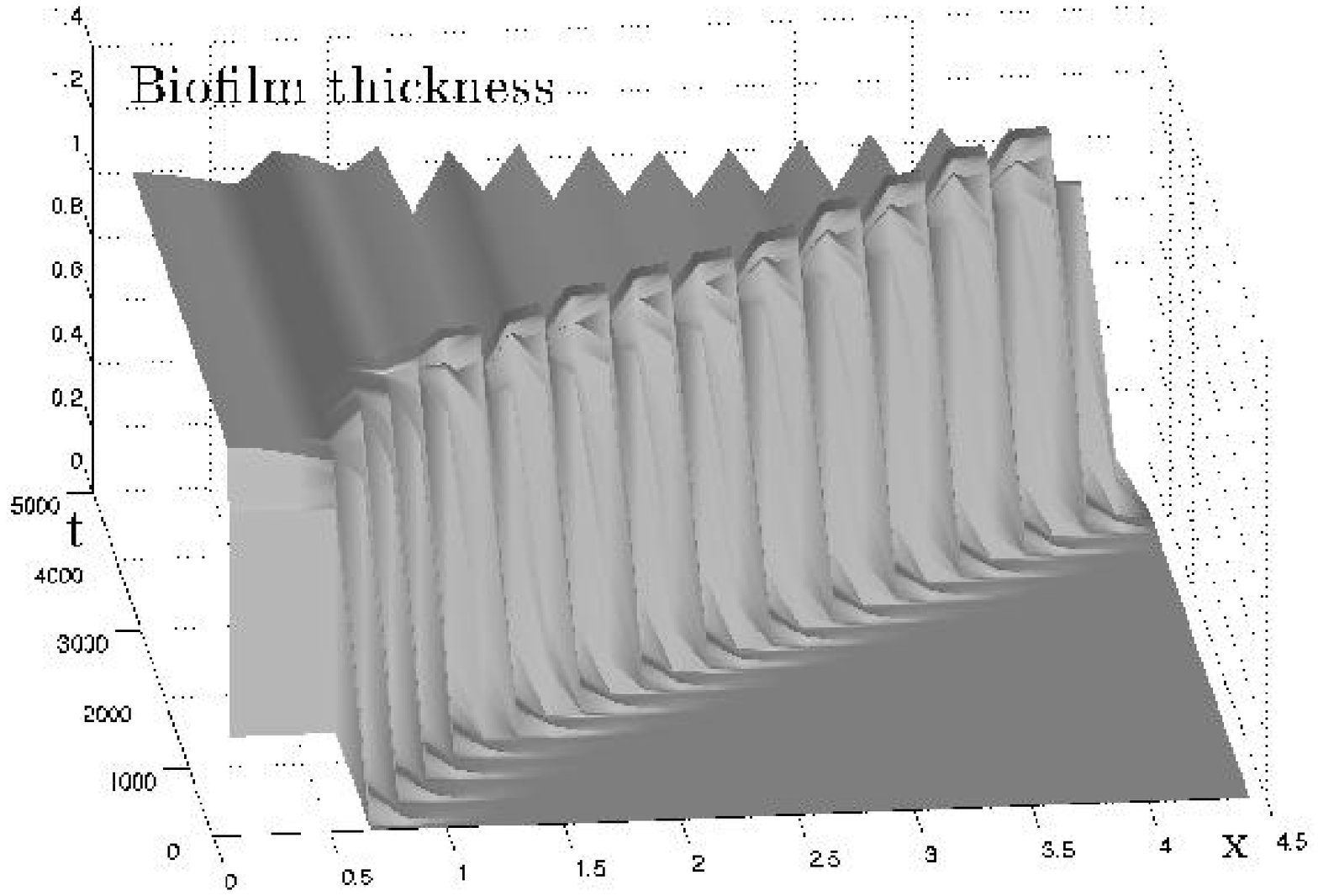}
\caption{\emph{Simulation 3 - Summarize of the swarm phenomenon. }}
\label{figSim2F}
\end{figure}

\begin{figure}
\centering 
\includegraphics[width=12cm,height=7cm]{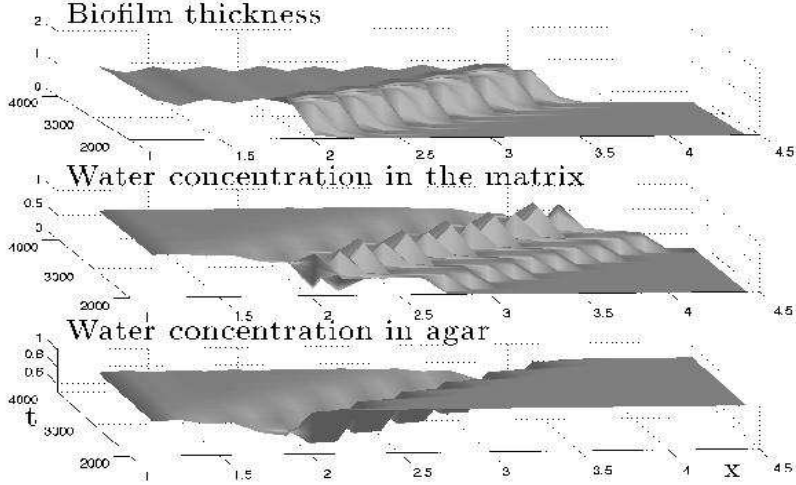}
\caption{\emph{Simulation 3 - Visualization of the relative positions of terraces, 
matrix water concentration maxima and agar water concentration minima.}}
\label{figSim2G}
\end{figure}
Figure \ref{figSim2B} shows the situation at the end of the colonization.
Here again, as in Figures \ref{figSim2F} and \ref{figSim2G} we may observe the 
relative regularity of the terraces.
\section{Conclusion}
\label{concl}
In this paper, we add a bio-physical slant in the \emph{Proteus mirabilis}
swarm description. 
The scientific method we use for this is as follows. 
We first express the bio-physical principles on which,
according to us, the swarm phenomenon could be based. 
Then, we translate those principles into a model
which is implemented. Simulations show that the model
behaviour may recover many of the aspects observed \emph{in situ}
and listed in the introduction.
This leads to the conclusion that the previously expressed bio-physical principles
are certainly the right ones to contribute to
a swarm explanation. Or, more modestly, those bio-physical 
principles are first way-marks on the path leading to a 
bio-physical understanding of the swarm phenomenon.

\:

\;

However, many questions need to be tackled.
From the modelling point of view, we notice that the 
model built here is dimensionless.
This means that variables have no physical unit 
and that fields are of magnitude order 1. 
Therefore, even if the model built here can qualitatively 
recover the swarm phenomenon, we are far from 
a model which can give quantitative results.
In order to build a model for quantitative simulations, 
we need to take into account the physical magnitudes of the 
phenomenon and to analyze physical constant ratios.

\:

From the mathematical point of view, it would be nice to
set out an existence result for the solution of the  
model built here.
Moreover, we may notice that this model
is instable. In particular, it has
been difficult to find values of the parameters that have led
to simulations showing terraces.
The stability properties, with respect to terrace generation or
to short-term evolution or to long-term behaviour
 of  the model
is certainly a very interesting topic which may bring out 
biological answers. 

\:

From the numerical point of view, the organic assumption concerning age and
time steps: $\bgest=\agest=\timest$, which is needed to make the numerical scheme
to be biomass preserving, will certainly result in bad computational results
if the present method is applied for operational simulations. Indeed, the simulation
of a real \emph{Proteus mirabilis} swarm experiment requires a time step which is much
smaller than the ageing characteristic time to resolve the dynamics in space.
Hence, it is important to find ways to reduce the induced cost.

The first way that can be considered, consists in noticing that the biomass preservation
property is always realised if the ageing process is done only every  $\prptab$
time steps, for an integer  $\prptab$, setting $\bgest=\agest= \prptab \timest$
and modifying the method in the following way. Initializing formula (\ref{A1}),
 (\ref{A3}), (\ref{A7}) and (\ref{A17AAA}) are left unchanged.
 \newline
Then for values of $n \in \prptab \nit^*$, formula (\ref{A2}) is applied removing term 
$\rrc_i^n$ and replacing $\timest$ by $\prptab \timest$ and  formula (\ref{A4})
 is applied replacing $\timest$ by $\prptab \timest$. Formula (\ref{A5}) to (\ref{A17000})
 and (\ref{A1600}) to (\ref{A17BBB}) are applied in state. 
\newline
For values of $n \notin \prptab \nit$, formula (\ref{A2}) is replaced by
\begin{gather}
\qc_i^{n+1} = \qc_i^n + \rrc_i^n \text{ if }n \in (\prptab \nit ^*+1) \; \text{ and } 
 \;  \qc_i^{n+1} = \qc_i^n \text{ otherwise },
\end{gather}
formula  (\ref{A4}) is not applied, formula (\ref{A5}) and (\ref{A6}) are
replaced by 
\begin{gather}
\zec^{n+1}_{k,i} = \zec^{n}_{k,i} \text{ for every } i \text{ and every } k,
\end{gather}
transfer condition (\ref{starstar}) is not applied,  (\ref{A8}) is replaced by
\begin{gather}
\label{A8modif}
\rhc^{n+1/2}_{k,p,i} = \rhc^{n}_{k,p,i},
\end{gather}
and (\ref{A13}) and  (\ref{A13.1}) are replaced by
\begin{gather}
\rhc^{n+1}_{k,p,i} = \rhc^{n+1/2}_{k,p,i} + \frac{\timest}{\posst}
(F^{n}_{k,p,i-1/2} - F^{n}_{k,p,i+1/2} ), \text{for every } p \text{ and } k,
\end{gather}
with fluxes given by  (\ref{A12}) and (\ref{A9}).
The definitions  (\ref{A15}) to (\ref{A152}) remain the same (but we can notice that
$\biomm^n_i$ do not evolve). Finally (\ref{A17000}), (\ref{A1600}), (\ref{A17}) 
and (\ref{A17BBB}) are applied in state.

The second way that needs to be investigated, consists in following the method
set out in Ayati \cite{Ayati2000} and which consists in deducing from equations
(\ref{S3}) and (\ref{S6}) evolution equations for functions $\zec(\tc,\age+\tc,\xc)$
and $\rhc(\tc,\age+t,\bge+t,\xc)$. From those equations we will be able to build age
discrete equations involving age discretization  possibly time dependant and
not requiring an organic link between time and age steps.

\:

From the biological point of view, it would be of great interest
to enrich the model in order to take into account genetic or
chemical events. In particular, it would be nice to understand
how gene expressions may regularize terrace generation. 

\bibliographystyle{plain}
\bibliography{biblio}

\end{document}